\documentclass[10pt]{article}
	\usepackage{caption}
	\usepackage{url}
	\usepackage{verbatim}
    \usepackage{graphicx}
    \usepackage{tikz}
    \usepackage{amsmath}
	\usepackage{mathrsfs}
	\usepackage{mathtools}
	\usepackage{stmaryrd}
    \usepackage{amssymb}
    \usepackage{amsthm}
    \usepackage{url}
    \usepackage{setspace}
    \usepackage{longtable}
    \usepackage{cmbright}
    \usepackage{color}
    \usepackage{fancyhdr}
    \usepackage{nicefrac}
	\usepackage{enumerate}
	\usetikzlibrary{patterns}
     \usepackage{indentfirst}
  \usepackage{hyperref}
  \hypersetup{
    colorlinks=true,
    linkcolor=blue,
    citecolor=blue,
}

    \usepackage{tabularx}

    \usepackage{multirow}
    \usepackage{floatrow}
    \usepackage{makecell}
    
    \usepackage{float}

    \newcommand{\bfx}{\boldsymbol{x}}

    \newcommand\addtag{\refstepcounter{equation}\tag{\theequation}}

    \numberwithin{equation}{section}

\usepackage{algorithm}
\usepackage{algorithmic}

\begin{document}

\title{Phase field modeling and computation of  vesicle growth or shrinkage}

\author{Xiaoxia Tang \thanks{Department of Applied Mathematics, Illinois Institute of Technology, Chicago, IL 60616 (xtang15@hawk.iit.edu)}
    \and
Shuwang Li\thanks{Department of Applied Mathematics, Illinois Institute of Technology, Chicago, IL 60616 (sli@math.iit.edu)}
	\and
John S. Lowengrub\thanks{Department of Mathematics, The University of California, Irvine, CA 92697 (jlowengr@uci.edu)}
	\and
Steven M. Wise\thanks{Corresponding author: Department of Mathematics, The University of Tennessee, Knoxville, TN 37996 (swise1@utk.edu)}}

	\maketitle
	\numberwithin{equation}{section}

	\begin{abstract}
We present a phase field model for vesicle growth or shrinkage induced by an osmotic pressure due to a chemical potential gradient. The model consists of an Allen-Cahn equation describing the evolution of phase field and a Cahn-Hilliard equation describing the evolution of concentration field. We establish control conditions for vesicle growth or shrinkage via a common tangent construction. During the membrane deformation, the model ensures total mass conservation and satisfies surface area constraint. We develop a nonlinear numerical scheme, a combination of nonlinear Gauss-Seidel relaxation operator and a V-cycles multigrid solver, for computing equilibrium shapes of a 2D vesicle.
Convergence tests confirm an $\mathcal{O}(t+h^2)$ accuracy.  Numerical results reveal that the diffuse interface model captures the main feature of dynamics: for a growing vesicle,  there exist circle-like equilibrium shapes if the concentration difference across the membrane and the initial osmotic pressure are large enough; while for a shrinking vesicle, there exists a rich collection of finger-like equilibrium morphologies.

	\end{abstract}


\section{Introduction}

Membranes considered in this paper are composed of bilayer lipid molecules with hydrophilic heads and two hydrophobic hydrocarbon chains. Lipid bilayers are the basic structural component of biological membranes. It is a semipermeable barrier to most solutes, including ions, proteins and other molecules. In an aqueous environment, a bilayer lipid membrane forms a vesicle (a closed bio-membrane containing fluid) to reduce the energy of the hydrophobic edges. Because of their relatively simple structure, vesicles are often used as a model system for studying fundamental physics underlying complicated biological systems such as cells and microcapsules. In addition, vesicles have also been used as building blocks to engineer artificial cells, e.g. biochemical microreactors operating in physiological environments \cite{Elani2015}. 

\begin{figure}[H]
\includegraphics[width=0.65\textwidth]{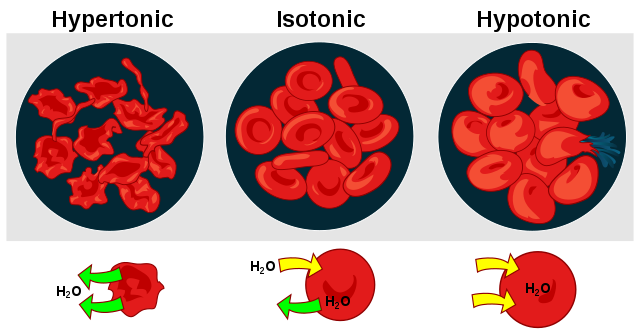}
\centering
\caption{Effect of different solutions on human red blood cells. The cell swells, shrinks or stays normal as water moves into or out of the cell down its concentration gradient in  hypotonic,  hypertonic, or isotonic solutions. From Wikipedia, the Free Encyclopedia. https://en.wikipedia.org/wiki/Tonicity \cite{wiki}.}
\label{fig:tonicity}
\end{figure}

Osmosis usually refers to the net movement of water molecules across a semipermeable membrane driven by a difference in concentration of solute on either side \cite{Albert2002, Strange2004}. Tonicity is another concept from osmosis.
It is operationally defined as the ability of a solution to shrink or swell specified cells. Hypotonicity describes any medium with a sufficiently low concentration of solutes to drive water to move into a cell due to osmosis. Hypertonicity describes any medium with a sufficiently high concentration of solutes to drive water to move out of a cell due to osmosis \cite{Albert2002,BAUMGARTEN2012}. In Figure \ref{fig:tonicity}, we show response of a human red blood cell to changes in tonicity of the extracellular fluid \cite{wiki}. Clearly, hypertonic solutions shrink cells, hypotonic solutions increase cell volume, and isotonic solutions neither swell nor shrink the cell. Motivated by these volume changes, in this paper, we develop a mathematical model to  simulate the effects of growth or shrinkage coming from the concentration gradient.

Mathematical modeling of membrane deformation has become an important area of research in biological and industrial system for a long time. At the continuum level, the mathematical description of vesicle conformation and deformation is a highly nonlinear, nonlocal moving boundary problem where the bilayer membrane serves as the moving boundary.
 Sharp interface models have been implemented to simulate the motion of vesicles in fluids \cite{VEERAPANENI2009A,VEERAPANENI2009B,SOHN2010119,Salac2011,Sohn2012,Shuwang2012,FRANK2013,Kai2016,Salac2018,Salac2022}.
For example, nonlinear wrinkling dynamics of a vesicle in an extensional flow and tumbling mechanism of two-dimensional vesicles in a shear flow are studied in \cite{Kai2014,Kai2017}.
Sharp interface models can satisfy the interface and inextensibility conditions exactly. That is,  the volume enclosed by membrane is automatically conserved for an incompressible fluid. The deformation based on the net mass transfer (gain or loss) across the membrane has not been considered. 
There are also sharp interface models of membrane deformation based on osmosis and diffusion \cite{Layton2006, Vogl2014, Jayathilake20101, Jayathilake20102, Mori2011, Yao2017,WANG2020,Quaife2021}. For example, boundary integral
simulations are used to investigate the effects of water permeability on the hydrodynamics of an inextensible membrane under a mechanical load in \cite{Quaife2021}.
An immersed boundary method \cite{PESKIN1977220}  for modeling convection and diffusion
of mass transfer through porous membranes under large deformations is proposed in a recent work \cite{WANG2020}. 

Phase field models have also been used to simulate the equilibrium configurations and the dynamics of vesicles \cite{DU2004450, Xiaoqiang2008,Lowengrub2009,Gu2016ATP}. It is well known that the phase field method is based on a diffusive interface approximation of a sharp interface, and treats the interface as a continuous but with steep change of properties of the two fluids by introducing a phase variable defined in the whole computational domain. That is, the interface is expressed by a thin internal transition layer of phase field \cite{Giga2017, Kobayashi2010}.
The phase field approach exhibits advantages in its simplicity in model formulation,  where the interface problem is posed as  a reaction-diffusion equation defined on the whole domain without requiring special treatment at the interface, thus ease of numerical implementation. Updating the interface position without explicitly tracking the interface is another attractive numerical feature of the phase field model \cite{Shen2012, Kobayashi2010, Provatas2010}.
The well known Allen-Cahn (AC) and Cahn-Hilliard (CH) equations are two gradient flow type PDEs describing the process of phase separation of a binary mixture, while the field variable is conserved in CH, but non-conserved in AC \cite{Giga2017, Bartels2015,Lee2014}. All spontaneous processes are accompanied by a decrease in free energy of the system.


In this paper, we develop a phase field model to simulate vesicle growth or shrinkage based on osmotic pressure, which arises due to a chemical potential gradient. In particular, we determine control conditions for growth and shrinkage via a common tangent construction. We simulate the growth and shrinkage effects subject to total mass conservation and surface area constraint, while allowing the mass exchange inside and outside the vesicle. 
Considering surface bending energy, osmotic pressure energy, and surface area constraint in addition to the surface energy used in the classical AC and CH equations, we derive an Allen-Cahn equation describing the evolution of phase field and a Cahn-Hilliard equation describing the evolution of concentration field.

The numerical computing and analysis are close to the method used in \cite{Wise2010}. We solve the problem by a nonlinear multigrid method \cite{Trottenberg2005, Henson2003, Kay2006}, which is a combination of nonlinear Gauss-Seidel relaxation operator and V-cycles multigrid solver. Specifically, at each time step we first need a nonlinear Gauss-Seidel smoothing operator, then use this smoothing operator on each hierarchical grid to get a better approximation. We demonstrate the  nearly optimal complexity of the multigrid solver and convergence of the scheme, which is of first order in time and second order in space. Numerical results reveal that for a growing vesicle,  there exist circle-like equilibrium shapes if the concentration difference across the membrane and the initial osmotic pressure are large enough; while for a shrinking vesicle, there exists a rich collection of finger-like equilibrium morphologies.

This paper is organized as follows. In Section \ref{sec: Model}, we define the model equations and analyse conditions for vesicle growth or shrinkage by common tangent construction. In Section \ref{sec:scheme}, the numerical scheme of the system and  nonlinear multigrid algorithm are presented. Numerical results  are given in Section \ref{sec:results}.

\section{Model Formulation}
\label{sec: Model}
\subsection{Evolution Equations}
We start by defining a Helmholtz free energy. Let $\Omega\subset\mathbb{R}^2$. The functions $\phi,\psi:\Omega\to \mathbb{R}$ are the phase fields (order parameters) describing the vesicle shape and the concentration of ionic fluid occupying the volume $\Omega$, respectively. $\left\{\mathbf{x}: \phi(\mathbf{x}) = 0\right\}$ determines the location of the membrane, while $\left\{\mathbf{x}: \phi(\mathbf{x}) = 1\right\}$ represents the interior phase (inside the vesicle),  and $\left\{\mathbf{x}: \phi(\mathbf{x}) = -1\right\}$ represents the exterior phase (outside the vesicle). We consider the following free energy densities \cite{Cahn1958,Giga2017, Du2005}:
	\begin{align}
f^{\rm surf}(\phi,\nabla\phi) &:= \frac{3\sqrt{2}}{4}\left(\frac{1}{\varepsilon}g(\phi)+\frac{\varepsilon}{2}|\nabla\phi|^2 \right) ,
	\\
f^{\rm bend}(\phi,\Delta\phi) & :=  \frac{3\sqrt{2}}{16\varepsilon} \left( \frac{1}{\varepsilon}g'(\phi) -\varepsilon\Delta\phi\right)^2 , 
	\\
f^{\rm osm}(\phi,\psi) & := \frac{1+p(\phi)}{2}f^{\rm in}(\psi) +\frac{1-p(\phi)}{2}f^{\rm out}(\psi),\label{eq:fosm}
	\end{align}
where $g$ is the standard  double-well function
$
g(\phi) = \frac{1}{4} \left(\phi^2 -1\right)^2
$. 
$\varepsilon$ is a small positive constant characterizing the thickness of the diffuse interface.
$f^{\rm in}(\psi)$ and $f^{\rm out}(\psi)$ are quadratic functions
	\[
f^{\rm in}(\psi) := \frac{\gamma_{\rm in}}{2} \left( \psi - \psi_{\rm in}\right)^2 + \beta_{\rm in} \quad \mbox{and} \quad f^{\rm out}(\psi) := \frac{\gamma_{\rm out}}{2} \left( \psi - \psi_{\rm out}\right)^2 + \beta_{\rm out}, \addtag \label{f_inout}
	\]
where $\gamma_{\rm in}, \psi_{\rm in}, \beta_{\rm in}, \gamma_{\rm out}, \psi_{\rm out}, \beta_{\rm out}$ are assumed to be positive parameters. $p$ is an interpolation function satisfying
	$p(1) = 1$ (interior phase) and $p(-1) = -1$ (exterior phase) and
	$
p'(-1) = p'(1) = 0
	$ 
as well. Now, we define free energies
	\begin{align}
F^{\rm surf} & := \int_\Omega \gamma_{\rm surf} f^{\rm surf}(\phi,\nabla\phi) \, d\bfx , 
	\\
F^{\rm bend} & := \int_\Omega \gamma_{\rm bend} f^{\rm bend}(\phi,\Delta\phi) \, d\bfx , 
	\\
F^{\rm area} & :=  \frac{\gamma_{\rm area}}{2} \left(\int_\Omega f^{\rm surf}(\phi,\nabla\phi)
 \, d\bfx - A  \right)^2 , 
 	\\
 F^{\rm osm} & := \int_\Omega f^{\rm osm} (\phi,\psi) \, d\bfx .
	\end{align}
Note that $F^{\rm surf}$ is related to the surface area of the vesicle, and $F^{\rm bend}$	is the surface bending energy. $A$ is the initial surface area of the vesicle, and $F^{\rm area}$ is the penalty term to numerically enforce the surface area constraint, since the surface area should remain unchanged for a vesicle with fixed amount of lipids. $\gamma_{\rm surf}$, $\gamma_{\rm bend}$, $\gamma_{\rm area}$ $> 0$. $F^{\rm osm}$ describes the osmotic energy arises in the mixture fluids with different concentrations.
The Helmholtz free energy is thus defined as
	\[
F[\phi,\psi] = F^{\rm surf}[\phi] + F^{\rm bend}[\phi] + F^{\rm area}[\phi] +  F^{\rm osm}[\phi,\psi]. \addtag
	\]

After energy variation, the dynamic equations read
\begin{align}
\partial_t \phi & = -M_\phi \mu, \label{eq:ACCHstart}
	\\
 \mu & = \delta_\phi F,
 	\\
 \partial_t \psi & = \nabla \cdot\left(M_\psi (\phi) \nabla \nu\right),
 	\\
 \nu & = \delta_\psi F,\label{eq:ACCHend}
 	\end{align}
where the phase variable $\phi$ is non-conserved and satisfies an Allen-Cahn equation, while the variable $\psi$ is conserved and satisfies a Cahn-Hilliard mass conservation equation. We take homogeneous Neumann boundary conditions on $\partial\Omega$. Here the mobility $M_\phi>0$ is a constant and $M_\psi(\phi) >0$ is a positive function of $\phi$. $\mu$ and $\nu$ are the chemical potentials
\begin{align}
	\begin{split}
	  \mu & = \gamma_{\rm surf} \frac{3\sqrt{2}}{4} \omega + \gamma_{\rm bend} \frac{3\sqrt{2}}{8}\left( \frac{\omega}{\varepsilon^2}g''(\phi) -  \Delta \omega\right) 
	\\
& \quad + \gamma_{\rm area} \left(\int_\Omega f^{\rm surf}(\phi,\nabla\phi)\, d\bfx -A \right) \frac{3\sqrt{2}}{4} \omega 
	\\
& \quad + \frac{p'(\phi)}{2} \left( f^{\rm in}(\psi) - f^{\rm out} (\psi)  \right) ,
	\end{split}\\
\omega & =  \frac{1}{\varepsilon} g'(\phi) -\varepsilon\Delta\phi  ,
	\\
\nu & = \frac{1+p(\phi)}{2}\cdot \frac{df^{\rm in}}{d\psi}(\psi) +\frac{1-p(\phi)}{2}\cdot\frac{df^{\rm out}}{d\psi}(\psi) .
	\end{align}
It is reasonable to assume that the mobility for $\psi$ degenerates in the interfacial region, since mass flow is limited to small channels in the vesicle membrane. Therefore, we assume that
	\[
M_\psi(\phi) = 1 - M_0(\phi^2-1)^2, \quad \mbox{for some} \ M_0 \in (0,1), \addtag
	\]
which follows that $M_\psi \ge 1-M_0 >0$. The system $\eqref{eq:ACCHstart}-\eqref{eq:ACCHend}$ is energy dissipative, and the dissipation rate is 
$$
\begin{aligned}
d_t F &=\int_{\Omega}\left\{\delta F_{\phi} \partial_t \phi+\delta F_{\psi} \partial_t \psi\right\} d \boldsymbol{x} \\
&=\int_{\Omega}\{\mu \partial_t \phi+\nu \partial_t \psi\} d \boldsymbol{x} \\
&=\int_{\Omega}\left\{\mu\left(-M_{\phi} \mu\right)+\nu\left(\nabla \cdot\left(M_{\psi} \nabla \nu\right)\right)\right\} d \boldsymbol{x} \\
&=\int_{\Omega}\left\{-M_{\phi}|\mu|^{2}-M_{\psi}|\nabla \nu|^{2}\right\} d \boldsymbol{x}\leq 0
\end{aligned}
$$

Next, we analyse the conditions for vesicle growth or shrinkage via the common tangent construction based on the osmotic free energy.  

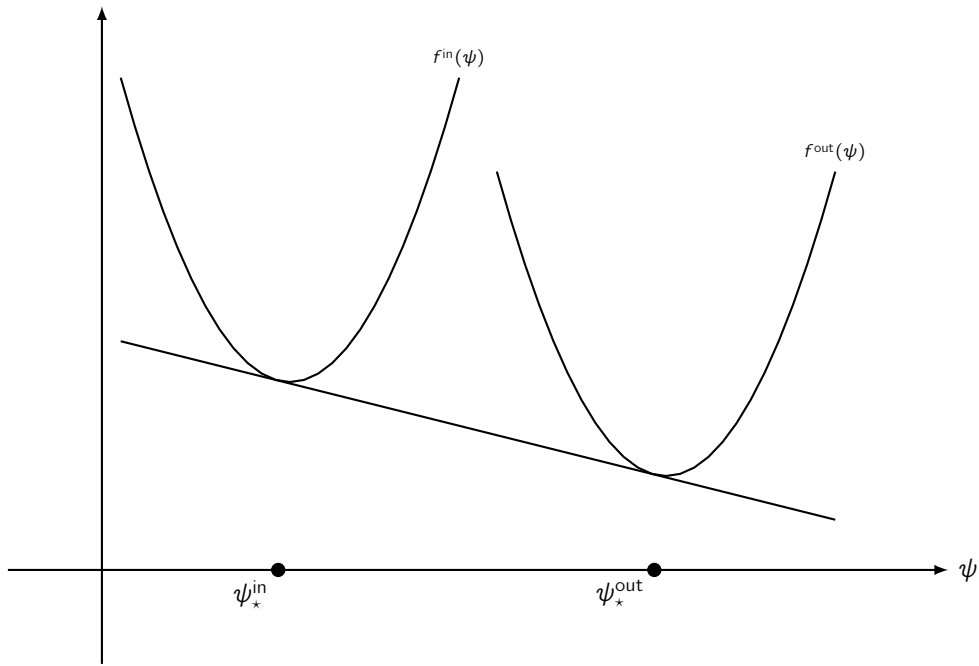
\begin{figure}[!htp]
	\begin{center}
	\begin{tikzpicture}[>=latex,scale=2.5]
    \draw[thick,->] (-0.5,0.0) -- (4.5,0.0) node[right] {$\psi$};
    \draw[thick,->] (0,-0.5) -- (0,3.0) node[above] {};
	\draw [thick,black,domain=0.1:1.9] plot (\x,{2.0*(\x-1.0)^2+1.0}) node[above] {$_{f^{\rm in}(\psi)}$};
	\draw [thick,black,domain=2.1:3.9] plot (\x,{2.0*(\x-3.0)^2+0.5}) node[above] {$_{f^{\rm out}(\psi)}$};
	\draw [thick,black,domain=0.1:3.9] plot (\x,{-0.25*\x+1/128+1+0.25*(1-1/16)}) node[above] { };
	\draw[black,fill=black] (1.0-1/16,0.0) node[below left] {$\psi^{\rm in}_\star$} circle (1pt);
	\draw[black,fill=black] (3.0-1/16,0.0) node[below left] {$\psi^{\rm out}_\star$} circle (1pt);
	\end{tikzpicture}
	\end{center}
\caption{Typical common tangent connection for the free energy densities for the interior $f^{\rm in}(\psi)$ (left) and exterior $f^{\rm out}(\psi)$ (right) \eqref{f_inout} phases. The vesicle will grow or shrink according to the initial (spatially uniform) states for $\psi$ inside and outside the vesicle. Suppose that $\psi^{\rm in}_{\star}$ and $\psi^{\rm out}_{\star}$ are the values of $\psi$ where the common tangent touches the respective free energy densities $f^{\rm in}(\psi)$ (left) and $f^{\rm out}(\psi)$ (right).}
	\label{fig:common-tangent}
	\end{figure}

	\subsection{Conditions for Growth or Shrinkage}
	\label{sec:3}
Suppose that $\psi^{\rm in}_{\star}$ and $\psi^{\rm out}_{\star}$ are the equilibrium concentration values for the interior and exterior regions obtained via the common tangent construction \cite{Provatas2010, Pelton2019}, and let us further assume that
	\[
0\le \psi^{\rm in}_{\star} < \psi^{\rm out}_{\star} \le 1, \addtag
	\]
as  shown in Figure~\ref{fig:common-tangent}. Now, suppose that we choose the following initial conditions for $\psi$:
	\[
\psi^{\rm in}(t=0) =: \psi^{\rm in}_0 < \psi^{\rm in}_{\star} \quad \mbox{and} \quad 	\psi^{\rm out}(t=0) =: \psi^{\rm out}_0 = \psi^{\rm out}_{\star}.  \addtag
	\]
This is the case that is illustrated in Figure~\ref{fig:growth}. The exterior phase is at its bulk equilibrium value, but the interior phase is not. The osmotic free energy is decreased as the concentration in the interior region goes up from the initial value $\psi^{\rm in}_0$ to the equilibrium  value $\psi^{\rm in}_{\star}$. Therefore mass will flow in until a global equilibrium is attained. In this case, the mass will be transferred into the  interior region (the volume of the vesicle will increase), though its surface area will remain unchanged.

	\begin{figure}[!htp]
	\begin{center}
	\begin{tikzpicture}[>=latex,scale=2.0,domain=0.2:3.2]
    \draw[thick,->] (-0.5,0.0) -- (4.5,0.0) node[right] {$x$};
    \draw[dashed,thick] (0.0,0.5) -- (2.0,0.5) node[right] {};
    \draw[dashed,red,ultra thick] (0.0,1.0) -- (2.666,1.0) node[right] {};
    \draw[dashed,red, ultra thick] (2.666,2.5) -- (4.5,2.5) node[right] {};
    \draw[dashed,thick] (2.0,2.5) -- (4.5,2.5) node[right] {};
    \draw[thick,->] (0,-0.5) -- (0,3.0) node[above] {};
    \draw[dashed,thick] (2.0,0.0)--(2.0,2.5) node[right] { };
    \draw[dashed,red,ultra thick] (2.666,0.0)--(2.6666,2.5) node[right] { };
	\draw[black,fill=black] (2.0,0.0) node[below left] {$x^{\rm int}_0$} circle (1pt);
	\draw[black,fill=black] (2.666,0.0) node[below left] {$x^{\rm int}_{\rm final}$} circle (1pt);
	\draw[black,fill=black] (0.0,0.5) node[below left] {$\psi^{\rm in}_0$} circle (1pt);
	\draw[black,fill=black] (0.0,1.0) node[below left] {$\psi^{\rm in}_\star$} circle (1pt);
	\draw[black,fill=black] (0.0,2.5) node[below left] {$\psi^{\rm out}_0=\psi^{\rm out}_\star$} circle (1pt);
	\fill[pattern=north west lines, pattern color=red] (0,0.5) rectangle (2,1.0);
	\fill[pattern=north west lines, pattern color=black] (2,1.0) rectangle (2.666,2.5);
	\node[text width=2cm] at (1.0,1.5){inner phase};
	\node[text width=2cm] at (4.0,1.5){outer phase};
	\end{tikzpicture}
	\end{center}
\caption{An example of volumetric growth. The black dashed regions show the initial state, and the red shows the final state. $x^{\rm int}_0$ is the initial position of the interface; $x^{\rm int}_{\rm final}$ is the position of the interface after global equilibrium is attained. The total  mass is conserved; mass from the black region is transferred into the red region.  The interface moves to the right, and the concentration in the inner phase increases.}	
	\label{fig:growth}
	\end{figure}

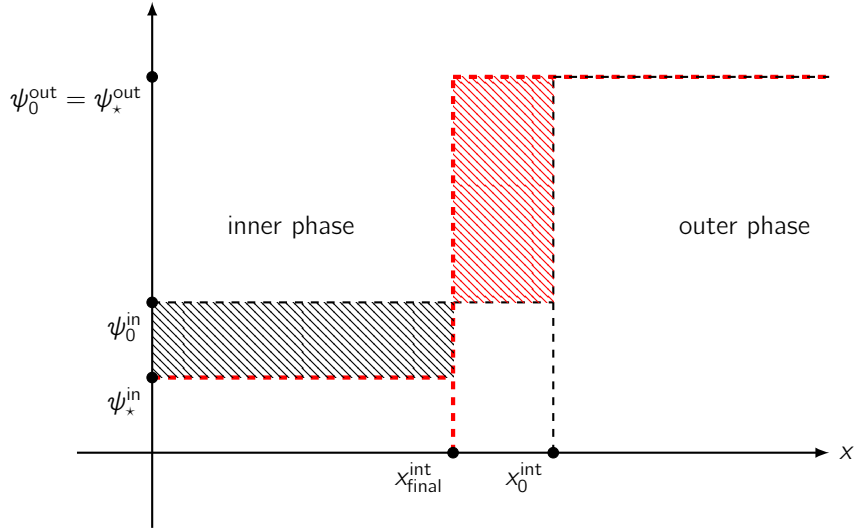
\begin{figure}[!htp]
	\begin{center}
	\begin{tikzpicture}[>=latex,scale=2.0,domain=0.2:3.2]
    \draw[thick,->] (-0.5,0.0) -- (4.5,0.0) node[right] {$x$};
    \draw[dashed,red,ultra thick] (0.0,0.5) -- (2.0,0.5) node[right] {};
    \draw[dashed,thick] (0.0,1.0) -- (2.666,1.0) node[right] {};
    \draw[dashed,thick] (2.0,2.5) -- (4.5,2.5) node[right] {};
    \draw[dashed,red, ultra thick] (2.0,2.5) -- (4.5,2.5) node[right] {};
    \draw[dashed,black, thick] (2.666,2.5) -- (4.5,2.5) node[right] {};
    \draw[thick,->] (0,-0.5) -- (0,3.0) node[above] {};
    \draw[dashed,red,ultra thick] (2.0,0.0)--(2.0,2.5) node[right] { };
    \draw[dashed,thick] (2.666,0.0)--(2.6666,2.5) node[right] { };
	\draw[black,fill=black] (2.0,0.0) node[below left] {$x^{\rm int}_{\rm final}$} circle (1pt);
	\draw[black,fill=black] (2.666,0.0) node[below left] {$x^{\rm int}_0$} circle (1pt);
	\draw[black,fill=black] (0.0,1.0) node[below left] {$\psi^{\rm in}_0$} circle (1pt);
	\draw[black,fill=black] (0.0,0.5) node[below left] {$\psi^{\rm in}_\star$} circle (1pt);
	\draw[black,fill=black] (0.0,2.5) node[below left] {$\psi^{\rm out}_0=\psi^{\rm out}_\star$} circle (1pt);
	\fill[pattern=north west lines, pattern color=black] (0,0.5) rectangle (2,1.0);
	\fill[pattern=north west lines, pattern color=red] (2,1.0) rectangle (2.666,2.5);
	\node[text width=2cm] at (1.0,1.5){inner phase};
	\node[text width=2cm] at (4.0,1.5){outer phase};
	\end{tikzpicture}
	\end{center}
\caption{An example of volumetric shrinkage. The black dashed regions show the initial state, and the red shows the final state. $x^{\rm int}_0$ is the initial position of the interface; $x^{\rm int}_{\rm final}$ is the position of the interface after global equilibrium is attained. The mass is conserved; mass from the black region is transferred into the red region. The interface moves to the left and the concentration in the inner phase decreases.}	
	\label{fig:shrinkage}
	\end{figure}

The shrinkage case is analogous and is illustrated in Figure~\ref{fig:shrinkage}. Suppose the initial conditions for this case are
	\[
\psi^{\rm in}(t=0) =: \psi^{\rm in}_0 > \psi^{\rm in}_{\star} \quad \mbox{and} \quad 	\psi^{\rm out}(t=0) =: \psi^{\rm out}_0 = \psi^{\rm out}_{\star}.  \addtag
	\]
In this case, the mass will be transferred into the exterior region, and the interior region will shrink (the volume of the vesicle decreases with constant surface area) as the osmotic free energy decreasing,  while the concentration in the interior region decreases from its initial value $\psi^{\rm in}_0$ to the equilibrium value $\psi^{\rm in}_{\star}$. 

	\section{Numerical Method}
	\label{sec:scheme}
	
In this section, we use a backward-time central-space method for discretization to get a semi-implicit numerical scheme, then solve the discrete system by a nonlinear Full Approximation Scheme (FAS) multigrid method,  which is a combination of nonlinear Gauss-Seidel relaxation operator and V-cycles multigrid solver.

	 \subsection{Discretization of Time}
	 
We propose the following time-discrete, space-continuous scheme of \eqref{eq:ACCHstart} - \eqref{eq:ACCHend}
\begin{align} 
\phi^{k+1}-\phi^{k} =&-s M_{\phi} \mu^{k+1}, \label{eq1} \\ 
\begin{split}
\mu^{k+1}=& \gamma_{1} \omega^{k+1}+\gamma_{2}\left(\frac{\omega^{k+1}}{\varepsilon^{2}} g^{\prime\prime}\left(\phi^{k}\right)-\Delta \omega^{k+1}\right) 
+\gamma_{3}\left(B^{k}-A\right) \omega^{k+1} \\
&+\frac{p^{\prime}\left(\phi^{k}\right)}{2}\left[f^{\rm in}\left(\psi^{k}\right)-f^{\rm out}\left(\psi^{k}\right)\right],
\end{split} \label{eq:Bk_A}\\
\omega^{k+1}=&\frac{1}{\varepsilon} g^{\prime}(\phi^{k+1})-\varepsilon \Delta \phi^{k+1},  
\\ 
\psi^{k+1}-\psi^{k} =&s \nabla \cdot\left(M_{\psi}(\phi ^k) \nabla \nu ^{k+1}\right), \\ 
 \nu ^{k+1}=&\frac{1+p\left(\phi^{k}\right)}{2} \frac{d f^{\text {in}}}{d \psi}\left(\psi^{k+1}\right)+\frac{1-p\left(\phi^{k}\right)}{2} \frac{d f^{\text {out}}}{d \psi}\left(\psi^{k+1}\right), \label{eq2}
 \end{align}
where $s$ is the time step,
$\gamma_{1}=\gamma_{\rm surf} \cdot \frac{3 \sqrt{2}}{4},$
$\gamma_{2}=\gamma_{\rm bend} \cdot \frac{3 \sqrt{2}}{8},$
$\gamma_{3}=\gamma_{\rm area} \frac{3 \sqrt{2}}{4},$
$
B^{k} =\int_{\Omega} f^{\rm surf}\left(\phi^{k}, \nabla \phi^{k}\right) d x,
$
and $\partial_n\phi^{k+1}=\partial_n\mu^{k+1}=\partial_n\omega^{k+1}=\partial_n\psi^{k+1}=\partial_n\nu^{k+1}=0 $ on $\partial \Omega$.

       \subsection{Discretization of Two-Dimensional Space}
       \subsubsection{Notations and Definitions}
       
Here we follow the notations and definitions of grid functions and difference operators used in \cite{Wise2010}. Consider $\Omega=(0,L_x) \times (0,L_y) \subset \mathbb{R}^2,$ with $L_x = m \cdot h,L_y = n \cdot h$, where $h>0$ is the spatial resolution and $m,n$ are positive integers. First, let's denote
\begin{align}
C_{m}&=\left\{\left(i-\frac{1}{2}\right) \cdot h | i=1, \ldots, m\right\},\\
C_{\bar{m}}&=\left\{\left(i-\frac{1}{2}\right) \cdot h | i=0, \ldots, m+1\right\},\\
E_{m}&=\{i \cdot h | i=0, \ldots, m\}.
\end{align}

\noindent $C_{m}$ and $C_{\bar{m}}$ are sets of $\emph{cell-centered points}$ of the interval $[0,L_x]$. The two points in $C_{\bar{m}}\backslash C_{m}$ are called $\emph{ghost points}$. The elements of $E_{m}$ are called $\emph{edge-centered points} $ of $[0,L_x]$. Analogously, $C_{n}$ and $C_{\bar{n}}$ contain the cell-centered points of $[0,L_y]$, and $E_{n}$ is a uniform partition of $[0,L_y]$ of size $n$. We will consider cell-centered points as domain of our discretized functions, and need the notations of edge-centered points in the definition of difference operators.  We define the function spaces
\begin{align}
\mathcal{C}_{m \times n} &=\left\{\phi: C_{m} \times C_{n} \rightarrow \mathbb{R}\right\}, & & \mathcal{C}_{\bar{m} \times \bar{n}}=\left\{\phi: C_{\bar{m}} \times C_{\bar{n}} \rightarrow \mathbb{R}\right\}, \\
\mathcal{C}_{\bar{m} \times n} &=\left\{\phi: C_{\bar{m}} \times C_{n} \rightarrow \mathbb{R}\right\}, & & \mathcal{C}_{m \times \bar{n}}=\left\{\phi: C_{m} \times C_{\bar{n}} \rightarrow \mathbb{R}\right\}, \\
\mathcal{E}_{m \times n}^{\mathrm{ew}} &=\left\{f: E_{m} \times C_{n} \rightarrow \mathbb{R}\right\}, & & \mathcal{E}_{m \times n}^{\mathrm{ns}}=\left\{f: C_{m} \times E_{n} \rightarrow \mathbb{R}\right\}.
\end{align}

The functions of $\mathcal{C}_{m \times n}$, $\mathcal{C}_{\bar{m} \times \bar{n}}$, $\mathcal{C}_{\bar{m} \times n} $, and $\mathcal{C}_{m \times \bar{n}}$ are called \emph{cell-centered functions}. In component form these functions are identified via $\phi_{i, j}:=\phi\left(x_{i}, y_{j}\right),$ where $x_{i}=\left(i-\frac{1}{2}\right) \cdot h, y_{j}=\left(j-\frac{1}{2}\right) \cdot h,$ and $i$ and $j$ are integers. The functions of $\mathcal{E}_{m \times n}^{\mathrm{ew}}$ and $\mathcal{E}_{m \times n}^{\mathrm{ns}}$ are called \emph{east-west edge-centered functions} and \emph{north-south edge-centered} functions, respectively. In component form east-west edge-centered functions are identified via $f_{i+\frac{1}{2}, j}:=$ $f\left(x_{i+\frac{1}{2}}, y_{j}\right),$ and north-south edge-centered functions are identified via $f_{i, j+\frac{1}{2}}:=f\left(x_{i}, y_{j+\frac{1}{2}}\right)$, where $x_{i+\frac{1}{2}}=i \cdot h, y_{j}=\left(j-\frac{1}{2}\right) \cdot h, x_{i}=\left(i-\frac{1}{2}\right) \cdot h, y_{j+\frac{1}{2}}=j \cdot h,$ and $i$ and $j$ are integers. Similarly, we define the edge-to-center difference operators $d_{x}: \mathcal{E}_{m \times n}^{\mathrm{ew}} \rightarrow \mathcal{C}_{m \times n}$ and $d_{y}: \mathcal{E}_{m \times n}^{\mathrm{ns}} \rightarrow$
$\mathcal{C}_{m \times n}$ component-wise via
\[
d_{x} f_{i, j}=\frac{1}{h}\left(f_{i+\frac{1}{2}, j}-f_{i-\frac{1}{2}, j}\right), \quad d_{y} f_{i, j}=\frac{1}{h}\left(f_{i, j+\frac{1}{2}}-f_{i, j-\frac{1}{2}}\right), \quad_{j=1, \ldots, n,}^{i=1, \ldots, m}  \addtag
\]
and center-to-center difference operators $c_{x}:\mathcal{C}_{\bar{m} \times n} \rightarrow \mathcal{C}_{m \times n}, c_{y}: \mathcal{C}_{m \times \bar{n}} \rightarrow \mathcal{C}_{m \times n}$,
\[
c_{x} \phi_{i, j}=\frac{1}{h}\left(\phi_{i+1, j}-\phi_{i, j}\right), \quad c_{y} \phi_{i, j}=\frac{1}{h}\left(\phi_{i, j+1}-\phi_{i, j}\right), \quad_{j=1, \ldots, n.}^{i=1, \ldots, m} \addtag
\]

\noindent The $x$ -dimension center-to-edge average and difference operators, respectively, $A_{x}, D_{x}:$ $\mathcal{C}_{\bar{m} \times n} \rightarrow \mathcal{E}_{m \times n}^{\mathrm{ew}}$ are defined component-wise as
\[
A_{x} \phi_{i+\frac{1}{2}, j}=\frac{1}{2}\left(\phi_{i, j}+\phi_{i+1, j}\right), \quad D_{x} \phi_{i+\frac{1}{2}, j}=\frac{1}{h}\left(\phi_{i+1, j}-\phi_{i, j}\right), \quad_{j=1, \ldots, n.}^{i=0, \ldots, m} \addtag
\]
Likewise, the $y$ -dimension center-to-edge average and difference operators, respectively,
$A_{y}, D_{y}: \mathcal{C}_{m \times \bar{n}} \rightarrow \mathcal{E}_{m \times n}^{\mathrm{ns}}$ are defined component-wise as
\[
A_{y} \phi_{i, j+\frac{1}{2}}=\frac{1}{2}\left(\phi_{i, j}+\phi_{i, j+1}\right), \quad D_{y} \phi_{i, j+\frac{1}{2}}=\frac{1}{h}\left(\phi_{i, j+1}-\phi_{i, j}\right), \quad_{j=0, \ldots, n.}^{i=1, \ldots, m} \addtag
\]
The standard 2D discrete Laplacian, $\Delta_{h}: \mathcal{C}_{\bar{m} \times \bar{n}} \rightarrow \mathcal{C}_{m \times n},$ is defined as
\begin{align}
\begin{split}
 \Delta_{h} \phi_{i, j} &=d_{x}\left(D_{x} \phi\right)_{i, j}+d_{y}\left(D_{y} \phi_{i, j}\right)\\
&=\frac{1}{h^{2}}\left(\phi_{i+1, j}+\phi_{i-1, j}+\phi_{i, j+1}+\phi_{i, j-1}-4 \phi_{i, j}\right), \quad_ {i=1, \ldots, m. 
}^{ j=1, \ldots, n}   
\end{split}
\end{align}
The numerical 2D integration, $B_h: \mathcal{C}_{m \times n}\rightarrow \mathbb{R}$ in \eqref{eq:Bk_A} is defined as
\begin{align}
\begin{split}
 B_h =& h^2 \sum_{i=1}^{m}\sum_{j=1}^{n} f^{\text{surf}}\left(\phi_{i,j}, \left(c_x\phi_{i,j},c_y\phi_{i,j}\right) \right)\\
=& h^2 \sum_{i=1}^{m}\sum_{j=1}^{n}
\frac{3\sqrt{2}}{4}\left\{\frac{1}{\varepsilon}g(\phi_{i,j})+\frac{\varepsilon}{2}\left[\left(\frac{\phi_{i+1, j}-\phi_{i, j}}{h}\right)^2 + \left(\frac{\phi_{i, j+1}-\phi_{i, j}}{h}\right)^2 \right] \right\}.
\end{split}
\end{align}
\noindent The numerical 2D surface area, $A_h: \mathcal{C}_{m \times n}\rightarrow \mathbb{R}$ in \eqref{eq:Bk_A} is defined as the numerical surface area of the initial data,
\[
A_h =B_h^0= h^2 \sum_{i=1}^{m}\sum_{j=1}^{n} f^{\text{surf}}\left(\phi_{i,j}^0, \left(c_x\phi_{i,j}^0,c_y\phi_{i,j}^0\right) \right). \addtag
\]

\subsubsection{Boundary Conditions}

In this paper, we use grid functions satisfying homogeneous Neumann boundary conditions on $\Omega$, that is, the cell-centered function $\phi \in \mathcal{C}_{\bar{m} \times \bar{n}} $ satisfies 
\begin{equation}
\begin{array}{l}
\phi_{0, j}=\phi_{1, j}, \quad \phi_{m+1, j}=\phi_{m, j}, \quad j=1, \ldots, n \\
\phi_{i, 0}=\phi_{i, 1}, \quad \phi_{i, n+1}=\phi_{i, n}, \quad i=0, \ldots, m+1,
\end{array}
\label{boundary conditions}
\end{equation}
\noindent we use the notation $\mathbf{n} \cdot \nabla_{h} \phi=0$ to indicate that $\phi$ satisfies \eqref{boundary conditions}.

\subsubsection{Fully-discrete Scheme}
With the notations defined above, the fully-discrete scheme for the equations \eqref{eq1}-\eqref{eq2} is: given $\phi^k, \psi^k \in \mathcal{C}_{\bar{m} \times \bar{n}}$, find the grid functions $\phi^{k+1},$ $\mu^{k+1},$ $\omega^{k+1},$ $\psi^{k+1},$ $\nu^{k+1}$ $\in \mathcal{C}_{\bar{m} \times \bar{n}}$ such that $\mathbf{n} \cdot \nabla_{h} \phi^{k+1}=\mathbf{n} \cdot \nabla_{h} \mu^{k+1}=\mathbf{n} \cdot \nabla_{h} \omega^{k+1}=\mathbf{n} \cdot \nabla_{h} \psi^{k+1}=\mathbf{n} \cdot \nabla_{h} \nu^{k+1}=0$, and
\begin{align}
    \phi^{k+1}-\phi^{k} =&-s M_{\phi} \mu^{k+1}, \label{eq3}\\ 
    \begin{split}
       \mu^{k+1}=&\gamma_{1} \omega^{k+1}+\gamma_{2}\left(\frac{\omega^{k+1}}{\varepsilon^{2}} g^{\prime\prime}\left(\phi^{k}\right)-\Delta_h \omega^{k+1}\right) 
+\gamma_{3}\left(B_h^{k}-A_h\right) \omega^{k+1}\\
&+\frac{p^{\prime}\left(\phi^{k}\right)}{2}\left[f^{\rm in}\left(\psi^{k}\right)-f^{\rm out}\left(\psi^{k}\right)\right], 
    \end{split}\\
\omega^{k+1}=&\frac{1}{\varepsilon} g^{\prime}(\phi^{k+1})-\varepsilon \Delta_h \phi^{k+1},\\ 
\psi^{k+1}-\psi^{k} =&s \left\{d_x\left(M_\psi \left(A_x\phi^k \right) D_x\nu^{k+1}\right)+ d_y\left(M_\psi \left(A_y\phi^k \right) D_y\nu^{k+1}\right)
\right\},\\ 
 \nu ^{k+1}=&\frac{1+p\left(\phi^{k}\right)}{2} \frac{d f^{\text {in}}}{d \psi}\left(\psi^{k+1}\right)+\frac{1-p\left(\phi^{k}\right)}{2} \frac{d f^{\text {out}}}{d \psi}\left(\psi^{k+1}\right).\label{eq4}
\end{align}

\subsection{Multigrid Solver}
\label{MG}
We now rewrite \eqref{eq3}-\eqref{eq4} as the following component form: find $\phi^{k+1}, \mu^{k+1}, \omega^{k+1}, \psi^{k+1},$ and $\nu^{k+1}$ in $\mathcal{C}_{\bar{m} \times \bar{n}}$ with boundary conditions $\mathbf{n} \cdot \nabla_{h} \phi^{k+1}=\mathbf{n} \cdot \nabla_{h} \mu^{k+1}=\mathbf{n} \cdot \nabla_{h} \omega^{k+1}=\mathbf{n} \cdot \nabla_{h} \psi^{k+1}=\mathbf{n} \cdot \nabla_{h} \nu^{k+1}=0$ whose components solve
\begin{align}
    \phi_{i,j}^{k+1}+s M_{\phi} \mu_{i,j}^{k+1} =&\phi_{i,j}^{k}, \label{eqApp1}\\ 
\mu_{i,j}^{k+1}-\left[\gamma_{1} +\frac{\gamma_{2}}{\varepsilon^{2}} g^{\prime\prime}\left(\phi_{i,j}^{k}\right)+\gamma_{3}\left(B_h^{k}-A_h\right)\right]\omega_{i,j}^{k+1}+\gamma_{2}\Delta_h \omega_{i,j}^{k+1}
=&\frac{p^{\prime}\left(\phi_{i,j}^{k}\right)}{2}\left[f^{\rm in}\left(\psi_{i,j}^{k}\right)-f^{\rm out}\left(\psi_{i,j}^{k}\right)\right],\\
\omega_{i,j}^{k+1}-\frac{1}{\varepsilon}\left(\left(\phi_{i,j}^{k+1}\right)^3-\phi_{i,j}^{k+1}\right)+\varepsilon \Delta_h \phi_{i,j}^{k+1} =&0, \\ 
\psi_{i,j}^{k+1}-s d_x\left(M_\psi \left(A_x\phi^k \right) D_x\nu^{k+1}\right)_{i,j}-s d_y\left(M_\psi \left(A_y\phi^k \right) D_y\nu^{k+1}\right)_{i,j}
 =&\psi_{i,j}^{k}, \\ 
 \begin{split}
 \nu_{i,j} ^{k+1}-\left[\frac{1+p\left(\phi_{i,j}^{k}\right)}{2}\cdot \gamma_{\rm in}+\frac{1-p\left(\phi_{i,j}^{k}\right)}{2}\cdot \gamma_{\rm out} \right]\psi^{k+1}_{i,j}=&-\frac{1+p\left(\phi_{i,j}^{k}\right)}{2} \gamma_{\rm in}\psi_{\rm in}\\
 & -\frac{1-p\left(\phi_{i,j}^{k}\right)}{2} \gamma_{\rm out}\psi_{\rm out}.
 \end{split} \label{eqApp2}
\end{align}

Let $\boldsymbol{\phi} = (\phi, \mu, \omega, \psi, \nu)^T$, define the $5\times m \times n$ nonlinear operator $\mathbf{N} = (N^{(1)},N^{(2)},N^{(3)},N^{(4)},N^{(5)})$ as
\begin{align}
    N_{i,j}^{(1)}&=\phi_{i,j}^{k+1}+s M_{\phi} \mu_{i,j}^{k+1},  \\ 
N_{i,j}^{(2)}&=\mu_{i,j}^{k+1}-\left[\gamma_{1} +\frac{\gamma_{2}}{\varepsilon^{2}} g^{\prime\prime}\left(\phi_{i,j}^{k}\right)+\gamma_{3}\left(B_h^{k}-A_h\right)\right]\omega_{i,j}^{k+1}+\gamma_{2}\Delta_h \omega_{i,j}^{k+1},
\\
N_{i,j}^{(3)}&=\omega_{i,j}^{k+1}-\frac{1}{\varepsilon}\left(\left(\phi_{i,j}^{k+1}\right)^3-\phi_{i,j}^{k+1}\right)+\varepsilon \Delta_h \phi_{i,j}^{k+1}, 
\\ 
N_{i,j}^{(4)}&=\psi_{i,j}^{k+1}-s d_x\left(M_\psi \left(A_x\phi^k \right) D_x\nu^{k+1}\right)_{i,j}-s d_y\left(M_\psi \left(A_y\phi^k \right) D_y\nu^{k+1}\right)_{i,j},
\\ 
N_{i,j}^{(5)}&= \nu_{i,j} ^{k+1}-\left[\frac{1+p\left(\phi_{i,j}^{k}\right)}{2}\cdot \gamma_{\rm in}+\frac{1-p\left(\phi_{i,j}^{k}\right)}{2}\cdot \gamma_{\rm out} \right]\psi^{k+1}_{i,j},
\end{align}
and the $5\times m \times n$ source $\mathbf{S} = (S^{(1)},S^{(2)},S^{(3)},S^{(4)},S^{(5)})$ as
\begin{align}
S_{i,j}^{(1)} &=\phi_{i,j}^{k}, \\
S_{i,j}^{(2)} &=\frac{p^{\prime}\left(\phi_{i,j}^{k}\right)}{2}\left[f^{\rm in}\left(\psi_{i,j}^{k}\right)-f^{\rm out}\left(\psi_{i,j}^{k}\right)\right],\\
S_{i,j}^{(3)} &=0, 
\\ 
S_{i,j}^{(4)} &=\psi_{i,j}^{k}, \\ 
S_{i,j}^{(5)} &=-\frac{1+p\left(\phi_{i,j}^{k}\right)}{2} \gamma_{\rm in}\psi_{\rm in}-\frac{1-p\left(\phi_{i,j}^{k}\right)}{2} \gamma_{\rm out}\psi_{\rm out}. 
\end{align}

\noindent Then, the system \eqref{eqApp1}-\eqref{eqApp2} is equivalent to $\mathbf{N}(\boldsymbol{\phi}^{k+1})=\mathbf{S}(\boldsymbol{\phi}^{k})$.

\par Next, we apply the nonlinear FAS multigrid method to solve $\mathbf{N}(\boldsymbol{\phi}^{k+1})=\mathbf{S}(\boldsymbol{\phi}^{k})$ for a given $\boldsymbol{\phi}^{k}$. The main points of this method are: (1) we first need a smoothing operator for generating smoothed approximate solutions of $\mathbf{N}(\boldsymbol{\phi})=\mathbf{S}$, here we use a nonlinear Gauss-Seidel relaxation operator; (2)  we then use this smoothing operator on each hierarchical grid to get better approximation of  $\boldsymbol{\phi}^{k+1}$. For further use, we represent the smoothing operator as
\[
\Bar{\boldsymbol{\phi}} = \text{Smooth}\left(\boldsymbol{\phi}, \mathbf{N}, \mathbf{S}, \lambda \right),\addtag
\]
where $\lambda$ is the number of smoothing sweeps.
Next, let's give the details of the relaxation. Here, $\ell$ is the index for iterative step, and we set
$$
\begin{array}{ll}
\phi_{i+\frac{1}{2}, j}^{\mathrm{ew}}:=A_{x} \phi_{i+\frac{1}{2}, j}^{k}, 
& \phi_{i, j+\frac{1}{2}}^{\mathrm{ns}}:=A_{y} \phi_{i, j+\frac{1}{2}}^{k} \\
{M_{\psi}}_{i+\frac{1}{2}, j}^{\mathrm{ew}}:=M_{\psi}\left(\phi_{i+\frac{1}{2}, j}^{\mathrm{ew}}\right), 
& {M_{\psi}}_{i, j+\frac{1}{2}}^{\mathrm{ns}}:={M_{\psi}}\left(\phi_{i, j+\frac{1}{2}}^{\mathrm{ns}}\right).\end{array}
$$
The Gauss-Seidel smoothing works as following: for every $(i,j)$, stepping lexicogrphically from $(1,1)$ to $(m,n)$, find $\phi_{i,j}^{\ell+1}, \mu_{i,j}^{\ell+1}, \omega_{i,j}^{\ell+1}, \psi_{i,j}^{\ell+1},$ and $\nu_{i,j}^{\ell+1}$ that solve

\begin{align}
    \phi_{i,j}^{\ell+1}+s M_{\phi} \mu_{i,j}^{\ell+1} =&S_{i,j}^{(1)}\left(\boldsymbol{\phi}^k\right), \label{eqApp3}\\
     \begin{split}
\mu_{i,j}^{\ell+1}-\left[\gamma_{1} +\frac{\gamma_{2}}{\varepsilon^{2}} g^{\prime\prime}\left(\phi_{i,j}^{k}\right)+\gamma_{3}\left(B_h^{k}-A_h\right)+\frac{4\gamma_2}{h^2}\right]\omega_{i,j}^{\ell+1}
=&S_{i,j}^{(2)}\left(\boldsymbol{\phi}^k\right)\\
&-\frac{\gamma_2}{h^2}\left[\omega_{i+1,j}^{\ell}+\omega_{i-1,j}^{\ell+1}+\omega_{i,j+1}^{\ell}+\omega_{i,j-1}^{\ell+1} \right],
\end{split}\\
\begin{split}
\omega_{i,j}^{\ell+1}-\left[\frac{1}{\varepsilon}(\phi_{i,j}^{\ell})^2+\frac{4\varepsilon}{h^2}\right]\phi_{i,j}^{\ell+1} =&S_{i,j}^{(3)}\left(\boldsymbol{\phi}^k\right) - \frac{1}{\varepsilon}\phi_{i,j}^{\ell}\\
&-\frac{\varepsilon}{h^2}\left[\phi_{i+1,j}^{\ell}+\phi_{i-1,j}^{\ell+1}+\phi_{i,j+1}^{\ell}+\phi_{i,j-1}^{\ell+1} \right],
\end{split} \label{eqApp3_2}
\\ 
\begin{split}
\psi_{i,j}^{\ell+1}+\frac{s}{h^2}\left[{M_{\psi}}_{i+\frac{1}{2}, j}^{\mathrm{ew}}+{M_{\psi}}_{i-\frac{1}{2}, j}^{\mathrm{ew}} +{M_{\psi}}_{i, j+\frac{1}{2}}^{\mathrm{ns}}+ {M_{\psi}}_{i, j-\frac{1}{2}}^{\mathrm{ns}}\right]\nu_{i,j}^{\ell+1}=& S_{i,j}^{(4)}\left(\boldsymbol{\phi}^k\right)+ \frac{s}{h^2}\left[{M_{\psi}}_{i+\frac{1}{2}, j}^{\mathrm{ew}} \nu_{i+1,j}^{\ell}+ {M_{\psi}}_{i-\frac{1}{2}, j}^{\mathrm{ew}} \nu_{i-1,j}^{\ell+1} \right.\\
&\left. + {M_{\psi}}_{i, j+\frac{1}{2}}^{\mathrm{ns}} \nu_{i,j+1}^{\ell}+  {M_{\psi}}_{i, j-\frac{1}{2}}^{\mathrm{ns}} \nu_{i,j-1}^{\ell+1}
\right],
\end{split} \label{eqApp4_1}
\\ 
 \nu_{i,j} ^{\ell+1}-\left[\frac{1+p\left(\phi_{i,j}^{k}\right)}{2}\cdot \gamma_{\rm in}+\frac{1-p\left(\phi_{i,j}^{k}\right)}{2}\cdot \gamma_{\rm out} \right]\psi_{i,j}^{\ell+1}=&S_{i,j}^{(5)}\left(\boldsymbol{\phi}^k\right). \label{eqApp4}
\end{align}
In practice, we use Cramer's Rule to solve this $5\times 5$ linear system \eqref{eqApp3}-\eqref{eqApp4} or use Cramer's Rule to solve $3\times 3$ linear system \eqref{eqApp3}-\eqref{eqApp3_2} and solve $2\times 2$ linear system \eqref{eqApp4_1}-\eqref{eqApp4} simultaneously, since \eqref{eqApp3}-\eqref{eqApp3_2} and \eqref{eqApp4_1}-\eqref{eqApp4} are independent.

Multigrid works on a hierarchy of grids. We use the smoothing operator on each level of the grids to get a  better approximation. Here, we set $\text{minlevel}\leq \text{level} \leq 0$, 0 means  the index of the finest grid, and minlevel is the index of the coarsest grid. We also need to transform the results between two levels of grids. By $\mathbf{I}_{\text{level}}^{\text{level-1}}$ we denote the restriction operator which is defined by cell-center averaging, and by $\mathbf{I}_{\text{level}-1}^{\text{level}}$ we denote the prolongation operator which is defined by piece-wise constant interpolation. $\mathbf{I}_{\text{level}}^{\text{level-1}}$ transfers fine grid functions to the coarse grid, while $\mathbf{I}_{\text{level}-1}^{\text{level}}$ transfers coarse grid functions to the fine grid. The following is the algorithm for our multigrid solver \cite{Wise2010}, in which 
$
\boldsymbol{\phi}_{\text{level}}^{k+1, m+1}=\text{FASVcycle}\left(\boldsymbol{\phi}_{\text{level}}^{k+1, m}, \mathbf{N}_{\text{level}}, \mathbf{S}_{\text{level}}, \lambda, \text{level} \right)
$ 
is the recursive FAS V-Cycle iteration operator and the superscript $m$ is the V-Cycle loop index.

\subsection{Algorithm}
\label{Algorithm}
\begin{center}
    \textbf{RECURSIVE FAS V-CYCLE OPERATOR}
\end{center}

$$
\boldsymbol{\phi}_{\text{level}}^{k+1, m+1}=\text{FASVcycle}\left(\boldsymbol{\phi}_{\text{level}}^{k+1, m}, \mathbf{N}_{\text{level}}, \mathbf{S}_{\text{level}}, \lambda, \text{level} \right)
$$

\noindent \textbf{Pre-smoothing:}
$$
\Bar{\boldsymbol{\phi}}_{\text{level}} = \text{Smooth}\left(\boldsymbol{\phi}_{\text{level}}^{k+1, m}, \mathbf{N}_{\text{level}}, \mathbf{S}_{\text{level}}, \lambda \right)
$$
\textbf{Coarse-grid correction:}

\par \textbf{If level$ >$ minlevel}
\begin{align}
&\mathbf{S}_{\text{level}-1} = \mathbf{I}_{\text{level}}^{\text{level}-1}\left(\mathbf{S}_{\text{level}}-\mathbf{N}_{\text{level}}\right)+ \mathbf{N}_{\text{level}-1} \left(\mathbf{I}_{\text{level}}^{\text{level}-1}  \Bar{\boldsymbol{\phi}}_{\text{level}}\right),\\
&\Bar{\boldsymbol{\phi}}_{\text{level}-1}=\text{FASVcycle}\left(\mathbf{I}_{\text{level}}^{\text{level}-1} \boldsymbol{\phi}_{\text{level}}, \mathbf{N}_{\text{level}-1}, \mathbf{S}_{\text{level}-1}, \lambda, \text{level}-1 \right),\\
&\hat{\boldsymbol{\phi}}_{\text{level}-1}=\Bar{\boldsymbol{\phi}}_{\text{level}-1}- \mathbf{I}_{\text{level}}^{\text{level}-1} \Bar{\boldsymbol{\phi}}_{\text{level}},\\
&\hat{\boldsymbol{\phi}}_{\text{level}}=\Bar{\boldsymbol{\phi}}_{\text{level}}+ I_{\text{level}-1}^{\text{level}} \hat{\boldsymbol{\phi}}_{\text{level}-1},\\
&\text{post-smooth: } \boldsymbol{\phi}_{\text{level}}^{k+1, m+1}=\text{Smooth}\left(\hat{\boldsymbol{\phi}}_{\text{level}}, \mathbf{N}_{\text{level}}, \mathbf{S}_{\text{level}}, \lambda \right).
\end{align}
\par \textbf{end if}

\noindent The combined algorithm of time stepping and the FAS V-Cycle iteration operator is given as follows. 

\begin{center}
    \textbf{COMBINED TIME STEPPING AND FAS V-CYCLE ITERATION ALGORITHM}
\end{center}
\textbf{Initialize} $\boldsymbol{\phi}_{0}^{k=0}$ \\[0.2cm]
Time Step Loop: \textbf{for} $k=0, k_{\max }-1$ \\[0.15cm]
\indent $\textbf{set } \boldsymbol{\phi}_{0 }^{k+1, m=0}=\boldsymbol{\phi}_{0}^{k}$\\[0.15cm]
\indent \textbf {calculate } $\mathbf{S}_{0}\left(\boldsymbol{\phi}_{0}^{k}\right)$
\\[0.15cm]
\indent V-cycle Loop: \textbf{for} $m=0, m_{\max }-1$
\\[0.15cm]
\indent \indent $\boldsymbol{\phi}_{0}^{k+1, m+1}=\text{FASVcycle}\left(\boldsymbol{\phi}_{0}^{k+1, m}, \mathbf{N}_{0}, \mathbf{S}_{0}, \lambda, 0 \right)$\\[0.15cm]
\indent \indent \textbf{if} $\left\|\mathbf{S}_{0}\left(\boldsymbol{\phi}_{0}^{k+1,m+1}\right)-\mathbf{N}_{0}\left(\boldsymbol{\phi}_{0}^{k+1, m+1}\right)\right\|_{2, \star}<\tau$, \textbf{then}\\[0.15cm]
\indent \indent \indent $\textbf{set } \boldsymbol{\phi}_{0}^{k+1}=\boldsymbol{\phi}_{0}^{k+1, m+1}$ and \textbf{exit} V-cycle Loop\\[0.15cm]
\indent \textbf{end for} V-cycle Loop\\[0.15cm]
\textbf{end for} Time Step Loop\\

\noindent Here $\tau>0$ is the stopping tolerance, and the norm is defined by 
\[
\|\mathbf{R}(\boldsymbol{\phi})\|_{2, *}:=\sqrt{\frac{1}{5mn} \sum_{k=1}^{5} \sum_{i=1}^{m} \sum_{j=1}^{n}\left(R_{i, j}^{(k)}(\boldsymbol{\phi})\right)^{2}}, \addtag \label{App:residual}
\]
where $\mathbf{R}(\boldsymbol{\phi}):=\mathbf{S}\left(\boldsymbol{\phi}^{k}\right)-\mathbf{N}(\boldsymbol{\phi})$ is the $5 \times m \times n$ residual array, and $R_{i, j}^{(k)}(\boldsymbol{\phi})$ are its components.

	\section{Numerical Results}
	\label{sec:results}
	
	In this section, we discuss numerical results of the discrete system \eqref{eq3}-\eqref{eq4} solved by the nonlinear FAS multigrid algorithm.  We present the results of convergence tests and perform sample computations. In all the tests below, we set  $L_x=L_y$ for simplicity, and use the interpolation function $p(\phi)=-\frac{1}{2}\phi^3+\frac{3}{2}\phi$ which satisfies 
	$p(1) = 1$, $p(-1) = -1$, and $p'(-1) = p'(1) = 0$ in \eqref{eq:fosm}.
	In the first set of tests in subsection  \ref{sec:5.1}, we show evidence that the multigrid solver converges with optimal (or near optimal) complexity. 
	In the second set of tests in subsection \ref{sec:5.2}, we provide evidence that the scheme is convergent and the global error is of first order in time and second order in space. We then present a series numerical studies on  the growth and shrinkage cases using the condition established in Section \ref{sec:3}.

\subsection{Convergence and Complexity of the Multigrid Solver}
\label{sec:5.1}
We perform six separate tests to demonstrate the convergence and near optimal complexity (with respect to the grid size $h$) of the multigrid solver. We provide evidence that the multigrid convergence rate is nearly independent of $h$. For all the tests we take the initial data
\[
    {\phi}^0_{i,j}= \tanh\left(\frac{0.18-\sqrt{0.75(x_i-0.5)^2+(y_j-0.5)^2}}{\sqrt{2}\varepsilon}\right) \label{ellipinitial1}, \addtag   
\]
\[
\psi^0_{i,j}=-\phi^0_{i,j} \times 0.1+0.7 \label{ellipinitial2}, \addtag
\]
and set the parameters $L_x= L_y=1.0, \gamma_{\rm surf}=1.0, \gamma_{\rm area}=1.0\times 10^{4},\gamma_{\rm in}=1.0\times 10^{5}, \gamma_{\rm out}=1.0\times 10^{5}, \psi_{\rm in}=0.1, \psi_{\rm out}=0.8, M_{0}=0.5, M_{\phi}=1.0, \beta_{\rm in}=0.0, \beta_{\rm out}=0.0.$  We use the temporal step size $s=5.0\times 10^{-7}$, and study the numerical results at the $20^{th}$ time step. We vary the spatial step size $h$ from $1.0/128$ to $1.0/1024$ and compare the number of multigrid iterations required to reduce the norm of the residual below the tolerance $\tau =1.0\times 10^{-8}$. Here, the stopping tolerance is $\|\mathbf{R}(\boldsymbol{\phi})\|_{2, \star} \leq \tau=1.0 \times 10^{-8}$, where $\mathbf{R}(\boldsymbol{\phi})$ and the norm are defined in \eqref{App:residual} of Section \ref{Algorithm}. $\lambda$ is the number of multigrid smoothing sweeps in the multigrid solver, as defined in Section \ref{MG}. Based on our experience as well as established in \cite{Trottenberg2005}, we expect that the optimal value of $\lambda$ should be less than 5. 

\begin{figure}[H]
\includegraphics[width=0.8\textwidth]{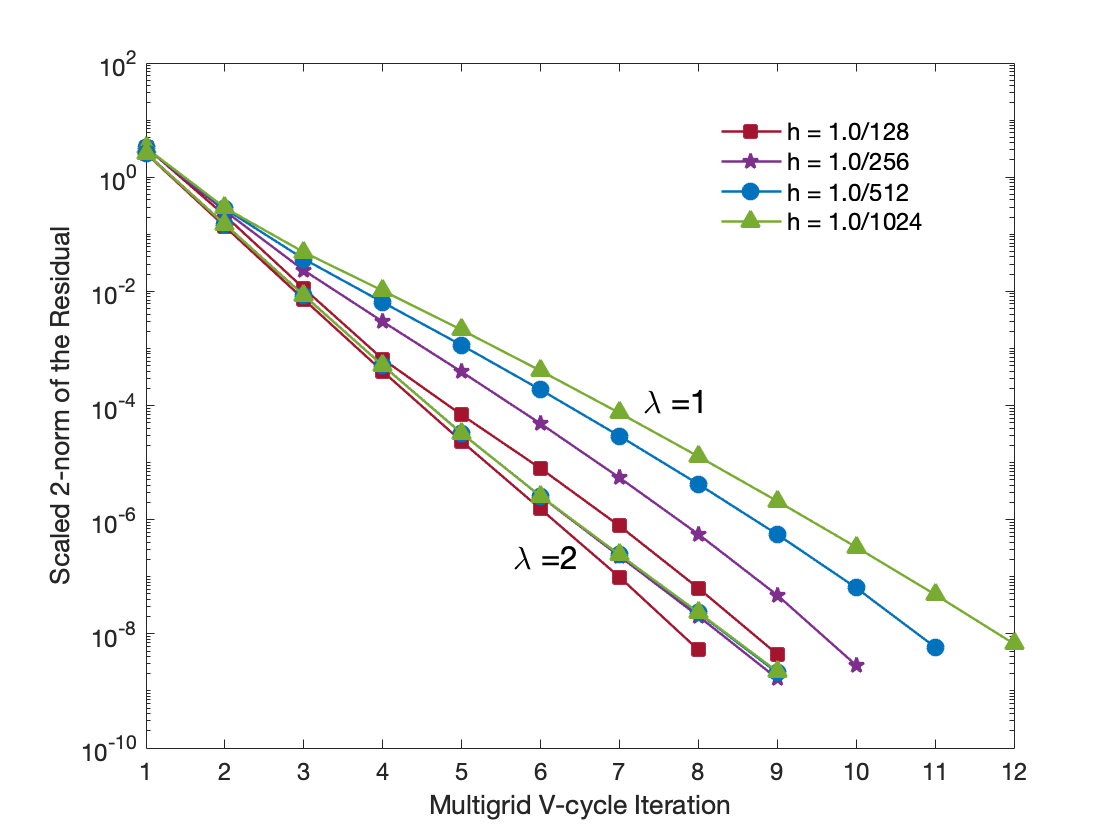}
\centering
\caption{The residual values per multigrid iteration in Test 2 and 6 at the $20^{th}$ time step with step size $s=5.0\times 10^{-7}$. The initial data are given by \eqref{ellipinitial1} - \eqref{ellipinitial2}.
The parameters are given in the text and in Table \ref{table:mg}. The results show that the residual reduction  is nearly independent of $h$ at $\lambda= 2$, which suggests the near optimal complexity of the solver}
\label{fig:mg}
\end{figure}

In Table~\ref{table:mg}, we show the number of multigrid iterations needed for various choices of $\varepsilon, \lambda$, and $\gamma_{\rm bend}$. We can see that for the smoothing parameters $\lambda= 2$, the required number of iterations is nearly independent of $h$. The detailed residual values for Test 2 and Test 6 in Table~\ref{table:mg} are given in graphical form in Figure \ref{fig:mg}, from which we observe for $\lambda= 2$ the norm of residual is reduced approximately the same factor at each iteration regardless of $h$. With $\lambda= 1$, we do not have this. By these features of multigrid operator with optimal complexity in \cite{Kay2006,Trottenberg2005}, it is evident that the multigrid solver here has near optimal complexity at $\lambda=2$.

\begin{table}[h!]
\centering
\begin{tabularx}{\textwidth}{>{\hsize=.2\hsize}XXXXXXXXXX}
\hline
                   & & Test 1 & Test 2 & Test 3  & Test 4  & Test 5 &Test 6   & Test 7 &Test 8                                     \\\hline
                   & $\varepsilon$ & $2\times 10^{-2}$ &$ 2\times 10^{-2}$ & $3\times 10^{-2}$ & $5\times 10^{-2}$ &$2\times 10^{-2}$ &$ 2\times 10^{-2}$ & $3\times 10^{-2}$ & $5\times 10^{-2}$ \\
                   & $\lambda $    & 1   & 1   &1 &1      & 2   &2   & 2 &2     \\ 
                   & $\gamma_{\rm bend}$ &0.1 &1 &0.5 &1 &0.1 &1 &0.5 &1       \\
                   \Xhline{1pt}
\multirow{4}{*}{$h$} 
                   & 1.0/128     & 11   & 9    & 8        & 11     & 8    & 8  & 10    & 7      \\
                   & 1.0/256     & 11   & 10   & 9        & 11     & 8    & 9  & 9     & 7       \\ 
                   & 1.0/512     & 10   & 11   & 10       & 12     & 9    & 9  & 9    & 8       \\ 
                   & 1.0/1024    & 11   & 12   & 11       & 12     & 9    & 9  & 9    & 8      \\
                   \hline                                      
\end{tabularx}
\caption{The number of multigrid iterations required to reduce the norm of the residual below the tolerance $\tau =1.0\times 10^{-8}$. The data are checked at the $20^{th}$ time step using the fixed temporal step size $s=5.0\times 10^{-7}$. The initial data are given by \eqref{ellipinitial1} - \eqref{ellipinitial2}.
The parameters are given in the text and in the table. The precise residual values for Test 2 and 6 are shown in Figure \ref{fig:mg}. Using the multigrid smoothing
parameter $\lambda = 2$ we observe, for a variety of parameter sets, that the required number of iterations is nearly independent of $h$} 
\label{table:mg}
\end{table}

       \subsection{Convergence of the Scheme as $\mathbf {s, h \rightarrow 0}$}
       \label{sec:5.2}
       
Next, we perform convergence tests of our scheme \eqref{eq3}-\eqref{eq4} as $s, h \rightarrow 0$. We expect that, at best, the global error in $\phi$ is $e_{t=T}=\mathcal{O}(s)+\mathcal{O}(h^2)$. To this end, we perform four tests similar to those in \cite{Wise2010}, under the same conditions except a refinement path of the form $s=Ch^2$. The initial data is given in \eqref{ellipinitial1}-\eqref{ellipinitial2} and the parameters used are $L_x= L_y=1.0, \varepsilon = 0.02, \gamma_{\rm bend}=0.1, \gamma_{\rm surf}=1.0, \gamma_{\rm area}=1.0\times 10^{4},\gamma_{\rm in}=1.0\times 10^{5}, \gamma_{\rm out}=1.0\times 10^{5}, \psi_{\rm in}=0.1, \psi_{\rm out}=0.8, M_{0}=0.5, M_{\phi}=1.0, \beta_{\rm in}=0.0, \beta_{\rm out}=0.0.$ and $T=6.4\times 10^{-4},$ where $T$ is the final time. We set time step size $s=6.4\times 10^{-6}$ and spatial step size $h =1/128,$ and check if the global error is reduced by a factor of 4 when $h$ is reduced by a factor of 2 and $s$ is reduced by a factor of 4.   Results in Table \ref{tabel:convergence} show evidence that the algorithm is convergent in space and the global error is indeed $e_{t=T}=\mathcal{O}(h^2)$. Refinement study in time step also show a first order accuracy.
In other
words, a global error of the form $e_{t=T}=\mathcal{O}(s)+\mathcal{O}(h^2)$
is consistent with the test results.

\begin{table}[h!]
\centering
 \begin{tabularx}{\textwidth}{X X X X X X }
 \hline
 Grid sizes  & $128^2-256^2$ &  & $256^2-512^2$& & $512^2-1024^2$  \\ [0.8ex] 
 
 Error &$1.0940\times 10^{-2}$ & & $2.8690\times 10^{-3}$ & &$7.2526\times 10^{-4}$\\ [0.8ex]
 Rate  & & 1.93 & & 1.98\\

 \hline
 \end{tabularx}
 \caption{Errors and convergence rates of the scheme \eqref{eq3}-\eqref{eq4}. Parameters are given in the text and the initial data is given by  \eqref{ellipinitial1}-\eqref{ellipinitial2}. Five tests with a refinement step size $s=Ch^2$ are presented. The error here is the global error between two nearby tests. The rate here suggests second order convergence rate in $h$ is attained, i.e., $e=\mathcal{O}(h^2)$} 
\label{tabel:convergence}
\end{table}
       \subsection{Vesicle Growth}
    In this section, we show the effect that the interior region of vesicle will expand, while the arclength (surface area in 3D) remain roughly a constant. We use initial conditions and parameters as the growth case described in Section \ref{sec:3}. Take initial condition $\phi^0_{i,j}$ as the smoothed result of $\hat{\phi}^0_{i,j}$ via a classical Cahn-Hilliard equation,
\begin{align}
 \hat{\phi}^0_{i,j}=&1,  \quad \quad if \quad (x_i,y_j)\in \Omega_1=\left\{(x_i,y_j)|(x_i-0.5)^2+(y_j-0.5)^2\leq r^2 \right\} \label{growthini_1}\\
 \hat{\phi}^0_{i,j}=&-1. \quad if \quad (x_i,y_j)\in \Omega \backslash \Omega_1
\end{align}
and
\[
  \psi^0_{i,j}=-\phi^0_{i,j} \times 0.35+0.45, \quad (x_i,y_j)\in \Omega \label{growthini_2} \addtag
\]
where 
$r = 0.18+0.03 \cos(10\theta), \theta \in [0, 2\pi]$. We use parameters $L_x= L_y=1.0, h=1.0/256, \varepsilon=0.01, \gamma_{\rm surf}=1.0, \gamma_{\rm bend}=0.05, \gamma_{\rm area}=5.0\times 10^{4},\gamma_{\rm in}=1.0\times 10^{5}, \gamma_{\rm out}=1.0\times 10^{5}, \psi_{\rm out}=0.8, M_{0}=0.5, M_{\phi}=1.0, \beta_{\rm in}=0.0, \beta_{\rm out}=0.0,$ $s=1\times 10^{-6},$ and the final time $T=2.5\times 10^{-2}.$ In this case, the equilibrium concentration values are $\psi_{\star}^{\rm in} = \psi_{\rm in}$ and $\psi_{\star}^{\rm out} = \psi_{\rm out}$.
The initial concentration of the outer phase is $\psi^0 = 0.8 = \psi_{\star}^{\rm out}$, which is at the equilibrium value. For the inner phase, the initial concentration is $\psi^0 = 0.1$. 
We next perform two sample computations with the equilibrium concentration of the inner phase  $\psi_{\star}^{\rm in} = 0.3$ or $\psi_{\star}^{\rm in} = 0.65$. According to the common tangent analysis in Section \ref{sec:3}, both cases will have growth effect.

\begin{figure}
\includegraphics[width=\textwidth]{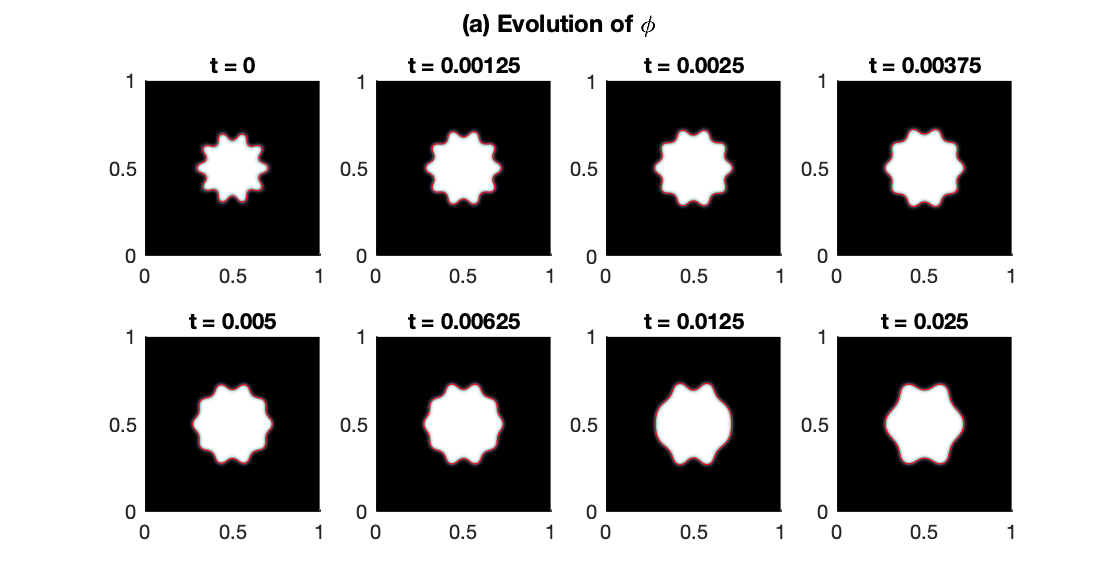}\\
\vspace{-0.2cm} 
\includegraphics[width=\textwidth]{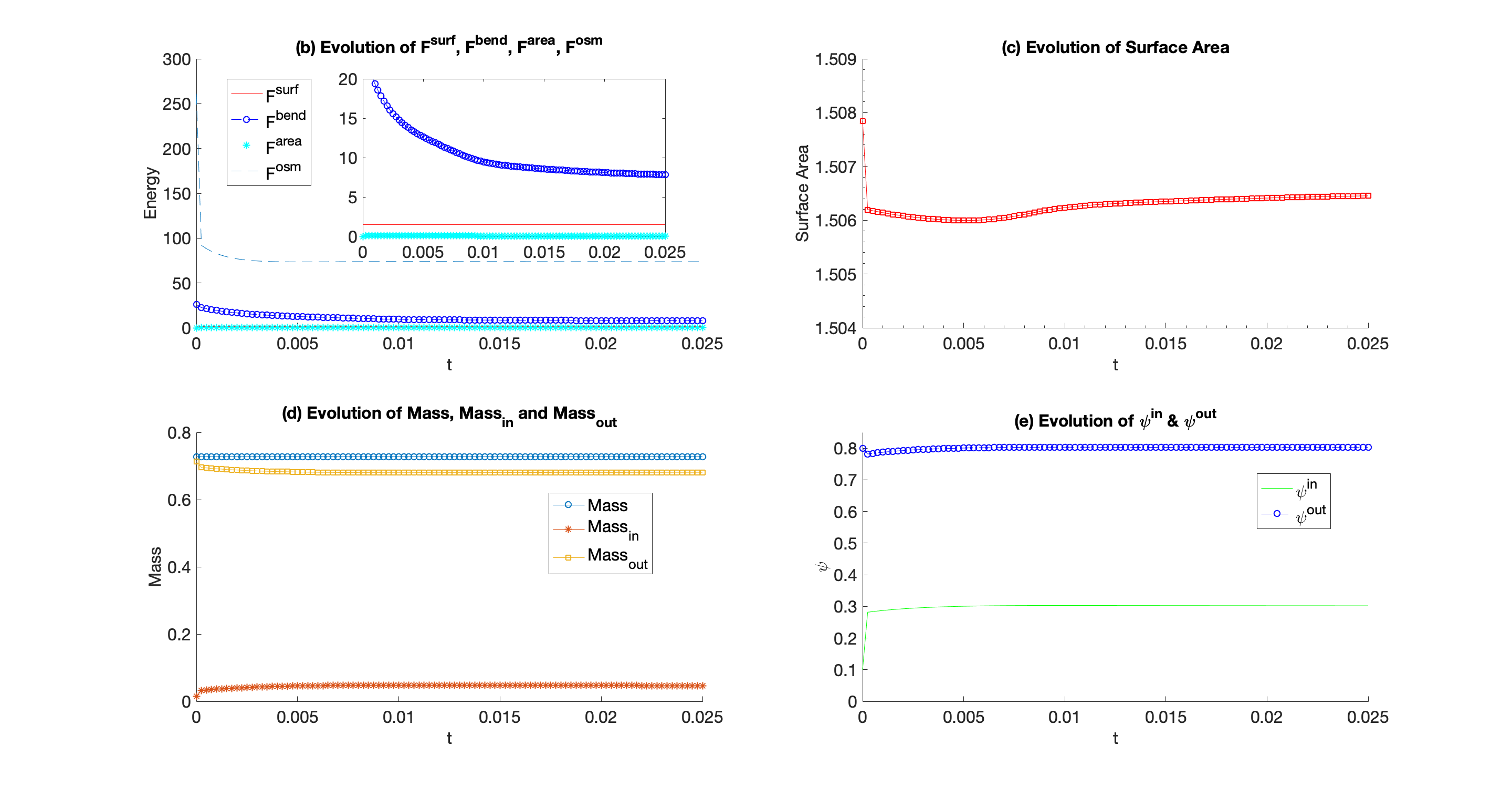}
\caption{The evolution of $\phi$, energy, surface area, mass, and concentration in a growth case. The initial condition is given in \eqref{growthini_1}-\eqref{growthini_2} and the parameters are in the text. The only initial difference with Figure \ref{fig:growth1} is that here we set $\psi_{\star}^{\rm in} = 0.3$. (a) shows the growth of inner region. The wrinkled interface stretches gradually with permanent area. In (b), the osmotic energy declines rapidly causing the interior growth. The bending energy drops during the swelling. The surface energy remain roughly a constant and $F^{\text{area}}$ stay close to 0 due to the penalty coefficient $\gamma_{\text{area}}$. (c) presents the change of $F^{\text{surf}}$, i.e., the interface area, is within $2\times 10^{-3}$, $0.133\%$ compared to the original surface area. In (d), the interior mass grows with the same amount lost in the exterior region, i.e., the total mass is conserved. (e) shows that $\psi^{\rm in}$ increases to the equilibrium value $\psi_{\star}^{\rm in}=0.3$, while $\psi^{\rm out}$ stays at the value of equilibrium $\psi_{\star}^{\rm out}=0.8$.}
\label{fig:growth_notcircle}
\end{figure}

\begin{figure}
\includegraphics[width=\textwidth]{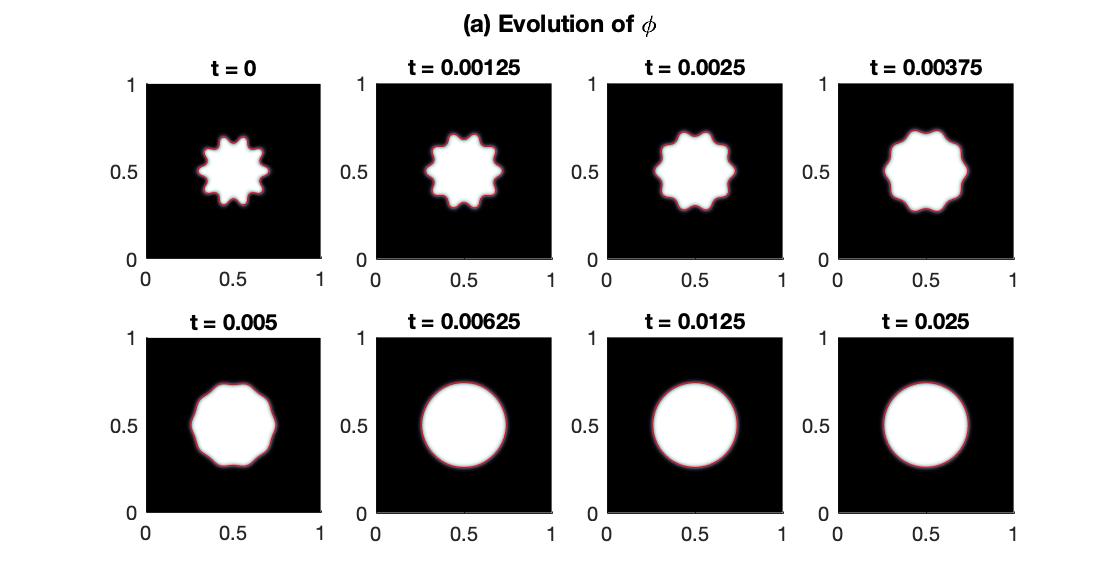}\\
\vspace{-0.2cm} \includegraphics[width=\textwidth]{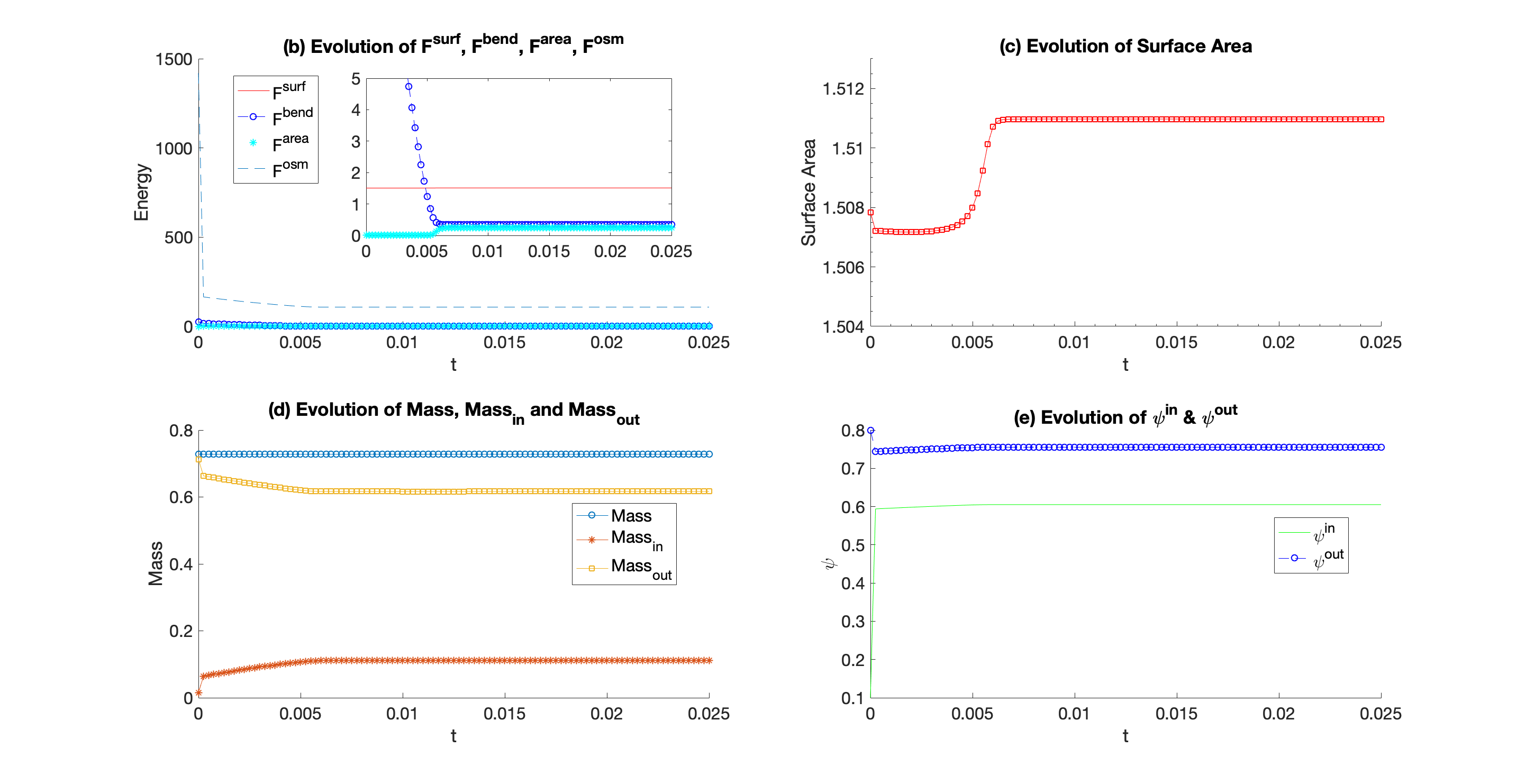}
\caption{The evolution of $\phi$, energy, surface area, mass, and concentration in a growth case. The initial condition is given in \eqref{growthini_1}-\eqref{growthini_2} and the parameters are in the text. The only initial difference with Figure \ref{fig:growth_notcircle} is that $\psi_{\star}^{\rm in} = 0.65$. (a) shows the growth of inner region. The wrinkled interface stretches gradually with permanent area and eventually the vesicle grows into a circle. In (b), the osmotic energy declines rapidly causing the interior growth. The bending energy drops until the wrinkled interface grows into a circle. The surface energy remain roughly a constant and $F^{\text{area}}$ stay close to 0 due to the penalty coefficient $\gamma_{\text{area}}$. (c) presents the change of $F^{\text{surf}}$, i.e., the interface area, is within $4\times 10^{-3}$, $0.265\%$ compared to the original surface area. In (d), the inside mass grows with the same amount lost in the outside region, i.e., the total mass is conserved. (e) shows that $\psi^{\rm in}$ increases toward the equilibrium value $\psi_{\star}^{\rm in}=0.65$ but stays at about 0.6 while $\psi^{\rm out}$ stays close to but not at the value of equilibrium $\psi_{\star}^{\rm out}=0.8$, because the vesicle already reaches a circular morphology (the shape with the maximum area with fixed arclength) and cannot accept any additional mass from the exterior region.}
\label{fig:growth1}
\end{figure}

In Figure \ref{fig:growth_notcircle} and Figure \ref{fig:growth1}, we show the result when $\psi_{\star}^{\rm in} = 0.3$ and $\psi_{\star}^{\rm in} = 0.65$, respectively.
In both calculations, we can see the growth of inner regions from the shape evolution (the evolution of $\phi$) in (a) of the two figures. The difference  is that the vesicle in Figure \ref{fig:growth1} grows into a circle, which is the state with maximized volume when the interface area stay unchanged; while the one in Figure \ref{fig:growth_notcircle} does not. Data in (b)-(e) of both figures explain how these changes and differences happen. 
(b) shows the evolution of energy $F^{\text{surf}}, F^{\text{bend}}, F^{\text{area}}$, and $F^{\text{osm}}$, respectively. It's easy to tell that the total energy drop significantly mainly due to the rapid decline of $F^{\text{osm}}$, while the surface energy $F^{\text{surf}}$ remain roughly constant because of the penalty coefficient $\gamma_{\text{area}}$. The detailed data of $F^{\text{surf}}$, i.e., the arclength of the interface, are shown in (c). We note that both changes are within $4\times 10^{-3}$, that is within $0.265\%$ compared to the original value. 
(d) shows the mass changes in the domain $\Omega$, the interior region (white), and the exterior  region (black), respectively. The interior mass grows with the same amount that is lost in the exterior region, i.e. the total mass is conserved. (e) shows the evolution of concentration of the interior region ($\psi^{\rm in}$) and the concentration of the exterior region ($\psi^{\rm out}$). Curves in Figure \ref{fig:growth_notcircle}(e) present that $\psi^{\rm in}$ increases to its equilibrium value $\psi_{\star}^{\rm in}=0.3$ and $\psi^{\rm out}$ stays at the equilibrium value $\psi^{\rm out}_{\star} = 0.8$, which is consistent with the analysis in Section \ref{sec:3}. But in Figure \ref{fig:growth1}, the outcomes are not exactly the same. The reason is that it needs certain amount of net mass moved  from the outer phase into the inner phase to attain the equilibrium value $\psi_{\star}^{\rm in} = 0.65$ in the inner phase, however the arclength inextensibility (interface area in 3D) here does not allow the inner concentration  to increase to the expected equilibrium value $0.65$,  because the vesicle already reaches a circular morphology, the shape with the maximum area (volume in 3D) for the fixed arclength, and cannot accept any additional mass from the exterior region.

\begin{figure}
\includegraphics[width=\textwidth]{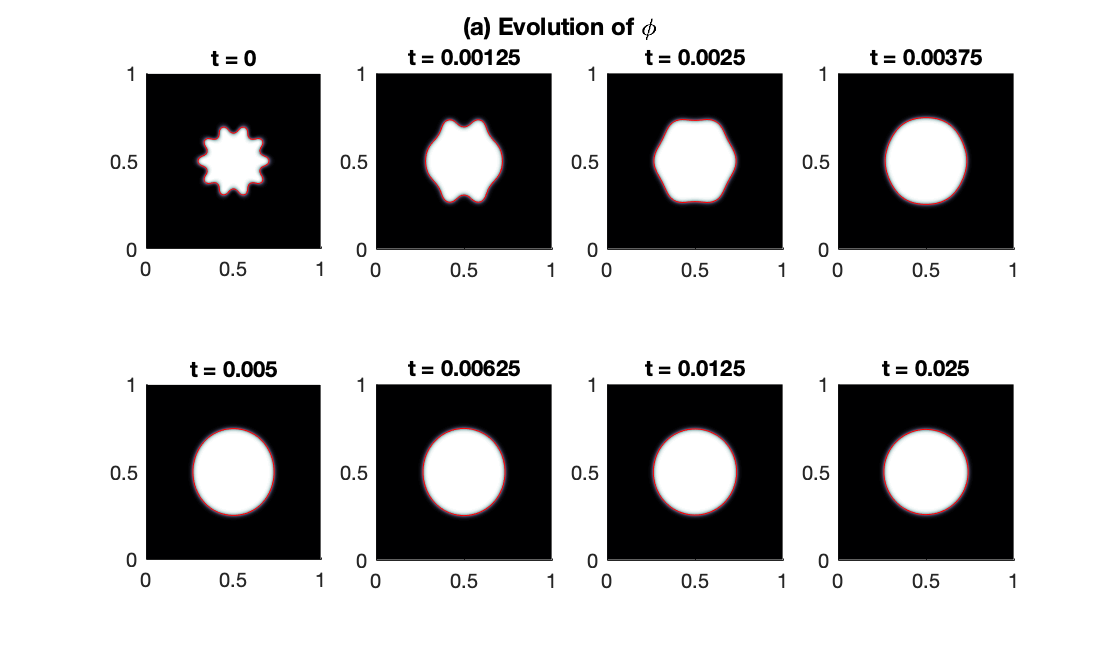}\\
\vspace{-0.2cm} 
\includegraphics[width=\textwidth]{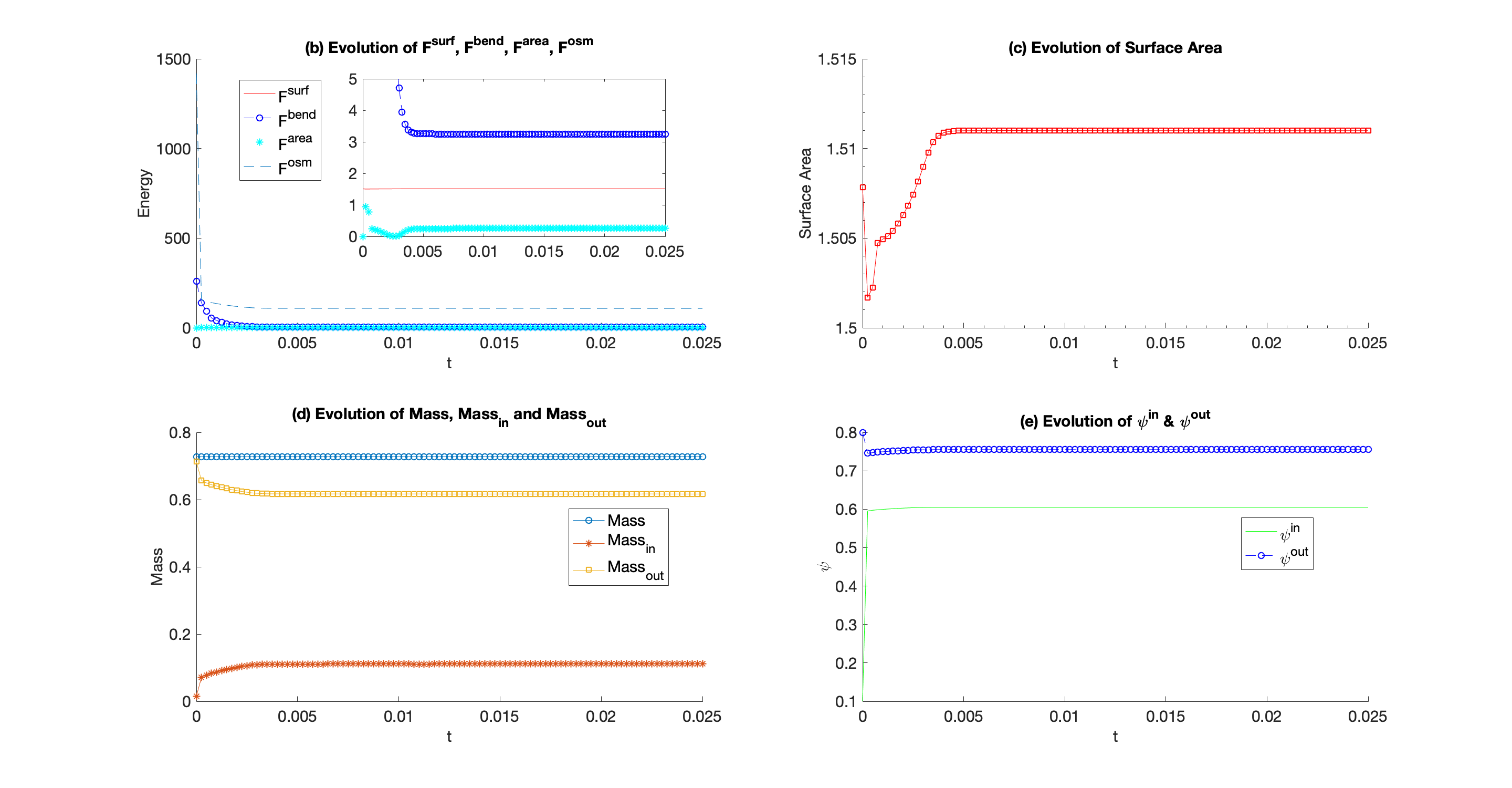}
\caption{The evolution of $\phi$, energy, surface area, mass, and concentration for another growth case. The initial condition and the parameters are the same as those in Figure \ref{fig:growth1} except $\gamma_{\rm bend}=0.5$. (a) shows the similar growth effect as in Figure \ref{fig:growth1}: the wrinkled interface stretches with permanent area and finally grows into a circle. But the shapes at the first several time steps are obviously different due to a faster decline of $F^{\text{bend}}$. In (b), the decline of osmotic energy causes the interior growth. The bending energy drops until the wrinkled interface grows into a circle. $F^{\text{surf}}$ remain roughly a constant and $F^{\text{area}}$ stay close to 0 due to the penalty coefficient $\gamma_{\text{area}}$. (c)-(d) present the surface area constraint and mass conservation. (e) shows $\psi^{\rm in}$ approaches the equilibrium value $\psi_{\star}^{\rm in}=0.65$, and $\psi^{\rm out}$ stay close to the value of the equilibrium $\psi_{\star}^{\rm out}=0.8$ as in Figure \ref{fig:growth1}.}
\label{fig:growth2}
\end{figure}

For the same initial condition in \eqref{growthini_1}-\eqref{growthini_2}, we now set $\gamma_{\rm bend}=0.5$ which means a ten times bending energy $F^{\text{bend}}$ compared to the previous calculation, and other parameters remain the same as those used in Figure \ref{fig:growth1}. The results are summarized  in Figure \ref{fig:growth2}. Compared to Figure \ref{fig:growth1}, the shape evolution in Figure \ref{fig:growth2} are obviously different at the first several time steps due to a much faster decline of $F^{\text{bend}}$. Eventually the interface evolves into a circle. 

\begin{figure}[H]
\includegraphics[width=0.8\textwidth]{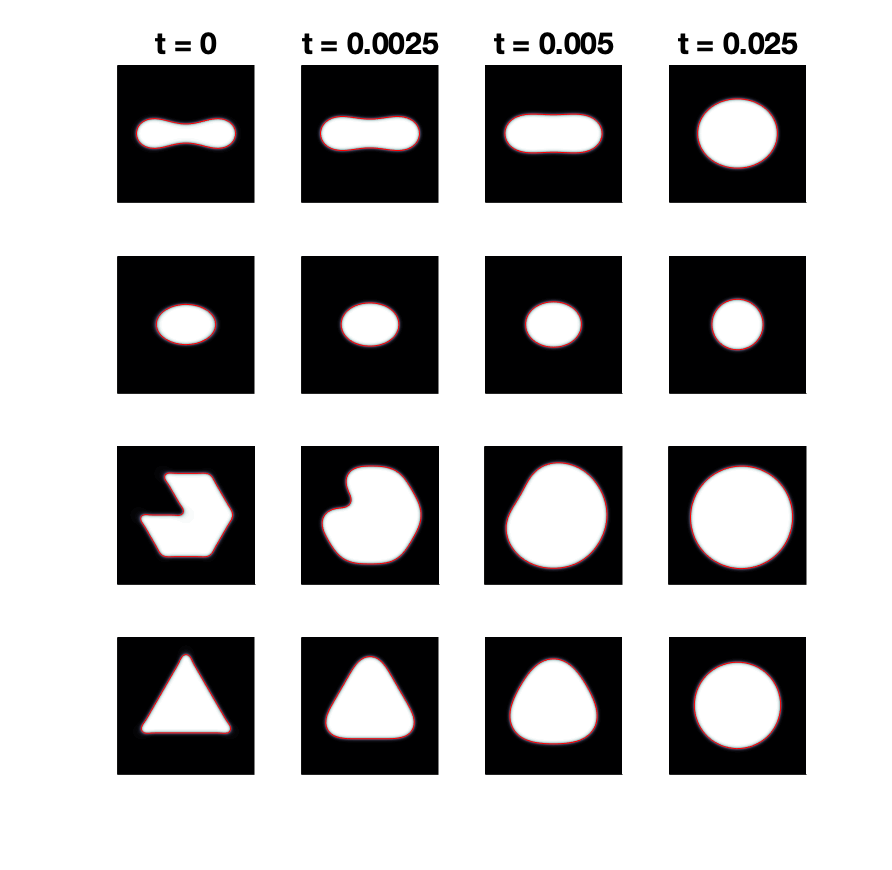}
\centering
\caption{Shape evolution of other growth examples with different initial conditions but same parameters as in Figure \ref{fig:growth1}. The white regions grow because of the decline of osmotic energy; the sharp corners swell faster for a fast drop of bending energy; as long as the inner equilibrium concentration $\psi_{\star}^{\rm in}$ and the initial osmotic energy $F^{\text{osm}}$ are great enough, the inner region will grow into a circle with preserved surface area.}
\label{fig:growth_other}
\end{figure}

In Figure \ref{fig:growth_other}, we show shape evolution of four vesicles with different initial configurations, and other parameters remain the same as those used in Figure \ref{fig:growth1}. Numerical experiments show that: the white region will grow with the decline of osmotic energy; the sharp corners will swell faster for a fast drop of bending energy; as long as the inner equilibrium concentration $\psi_{\star}^{\rm in}$ and the initial osmotic energy $F^{\text{osm}}$ are large enough, the inner region will grow into a circle with preserved arclength. 

\subsection{Shrinkage Model}
       
In this section,  we study the shrinking effect following the condition discussed in section \ref{sec:3}. That is the  area (volume in 3D) of the  vesicle will decrease with prescribed arclength (surface area in 3D). We take the initial shape  $\phi^0_{i,j}$ as a smoothed result of $\hat{\phi}^0_{i,j}$ via a classical Cahn-Hilliard equation,
\begin{align}
 \hat{\phi}^0_{i,j}=&1,  \quad \quad if \quad (x_i,y_j)\in \Omega_1=\left\{(x_i,y_j)|(x_i-0.5)^2+(y_j-0.5)^2\leq r^2 \right\} \label{shrinkini_1}\\
 \hat{\phi}^0_{i,j}=&-1. \quad if \quad (x_i,y_j)\in \Omega \backslash \Omega_1
\end{align}
and
\[
  \psi^0_{i,j}=-\phi^0_{i,j} \times 0.1+0.7, \label{shrinkini_2} \quad (x_i,y_j)\in \Omega \addtag
\]
where 
$r = 0.3+0.01 \cos(10\theta), \theta \in [0, 2\pi]$. We set parameters $L_x= L_y=1.0,h=1.0/256, \varepsilon=0.01, \gamma_{\rm surf}=1.0, \gamma_{\rm bend}=0.1, \gamma_{\rm area}=5.0\times 10^{4},\gamma_{\rm in}=1.0\times 10^{5}, \gamma_{\rm out}=1.0\times 10^{5}, \psi_{\rm in}=0.1, \psi_{\rm out}=0.8, M_{0}=0.5, M_{\phi}=1.0, \beta_{\rm in}=0.0, \beta_{\rm out}=0.0,$ $s=1.0\times 10^{-6}$, and the final calculation $T=4.0\times 10^{-2}.$ In this case, the equilibrium concentration values are $\psi_{\star}^{\rm in} = \psi_{\rm in}=0.1$ and $\psi_{\star}^{\rm out} = \psi_{\rm out}=0.8$.

\begin{figure}
\includegraphics[width=\textwidth]{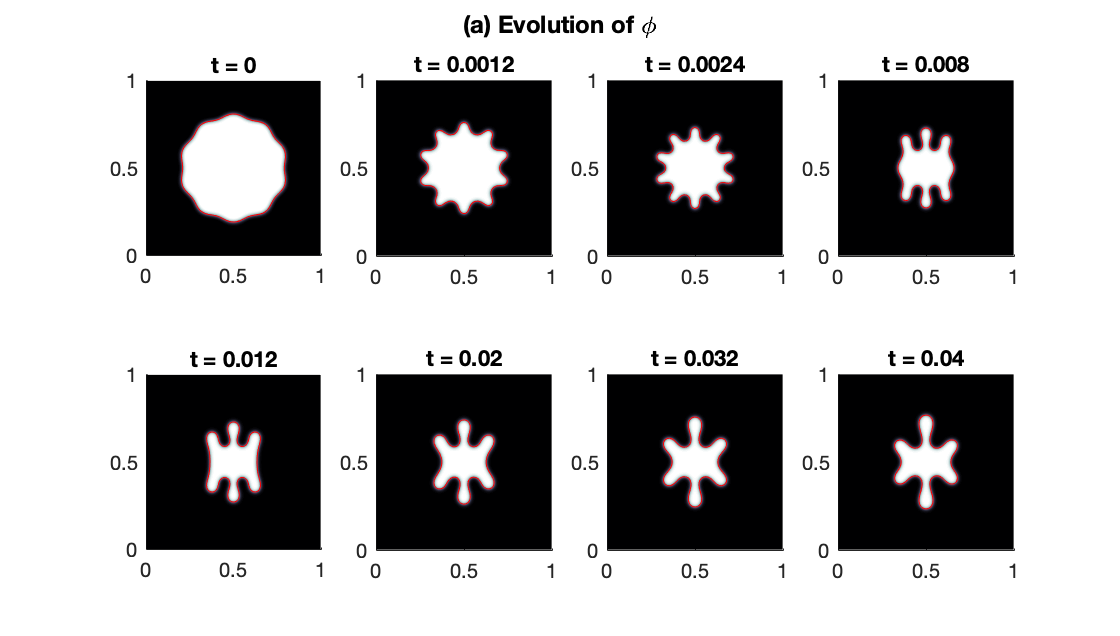}\\
\vspace{-0.2cm} 
\includegraphics[width=\textwidth]{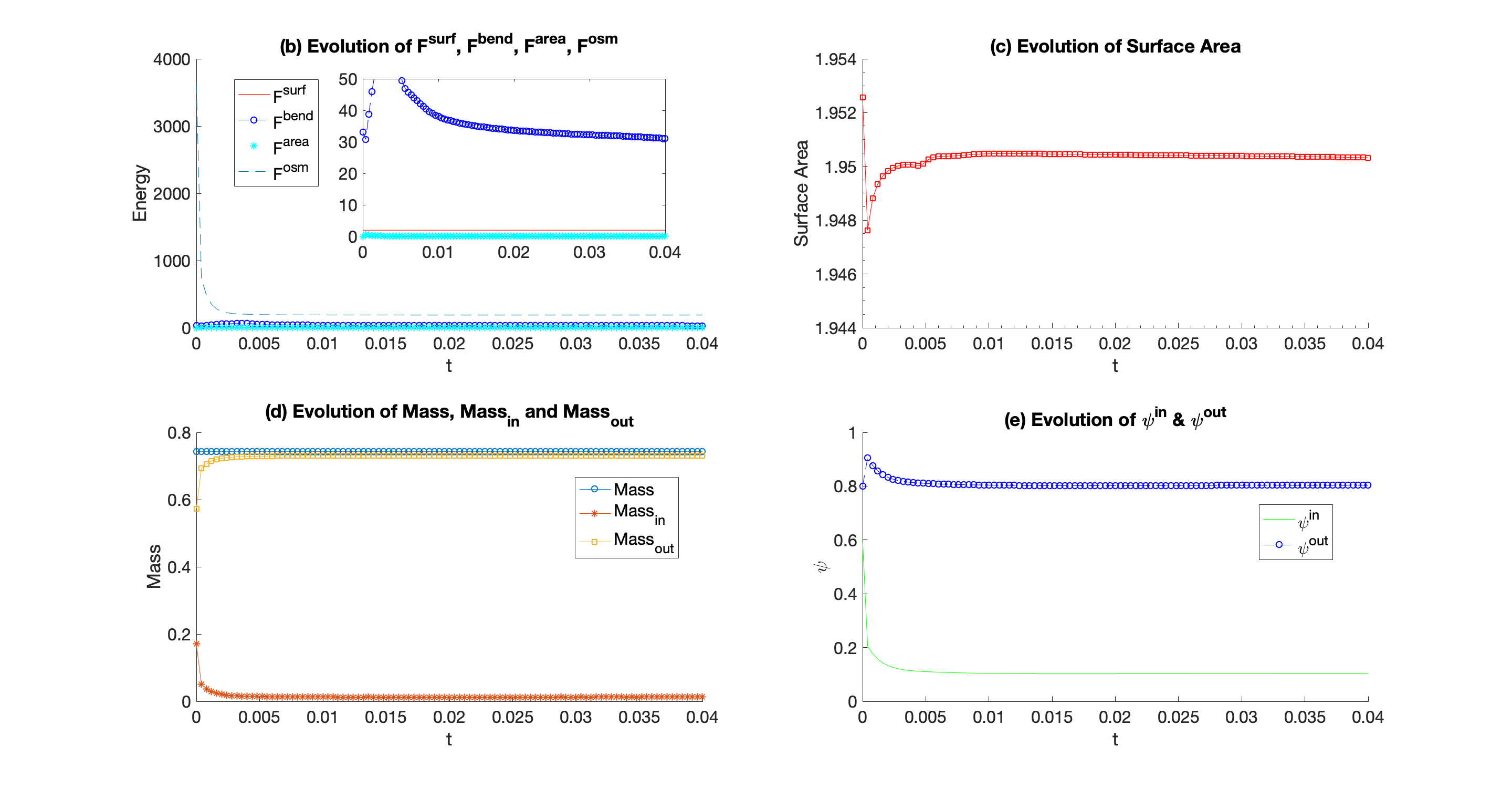}
\caption{The evolution of $\phi$, energy, surface area, mass, and concentration in a shrinkage case. The initial condition is given in \eqref{shrinkini_1}-\eqref{shrinkini_2} and the parameters are in the text. (a) shows shrinkage effects where the interface becomes wrinkled with prescribed surface area and the sharp corners shrink to form finger-like structures. In (b), the osmotic energy declines rapidly causing the inner region shrinking. The bending energy slightly grows at first because of the shrinking. The surface energy remain roughly a constant and $F^{\text{area}}$ stay close to 0 due to the surface area constraint. (c) presents the change of $F^{\text{surf}}$ is within $5\times 10^{-3}$, $0.26\%$ compared to the original surface area. (d) implies the mass conservation. In (e), $\psi^{\rm in}$ approaches the equilibrium value $\psi_{\star}^{\rm in}=0.1$, while $\psi^{\rm out}$ stays at the equilibrium value $\psi_{\star}^{\rm out}=0.8$.}
\label{fig:shrink1}
\end{figure}

In Figure \ref{fig:shrink1}(a), we show the shape evolution of $\phi$, which is obviously shrinking of the inner region while the interface becomes wrinkled. In Figure \ref{fig:shrink1}(b),  we show  the energy evolution curves of $F^{\text{surf}}, F^{\text{bend}}, F^{\text{area}}$, and $F^{\text{osm}}$. $F^{\text{bend}}$ slightly increases at early times when $F^{\text{osm}}$ drop significantly, leading to a shrinkage of the vesicle volume. $F^{\text{surf}}$, however,  stays roughly unchanged due to the arclength constraint. The detailed data of the surface area could be checked in Figure \ref{fig:shrink1}(c). Figure \ref{fig:shrink1}(d) gives the change of mass of the interior and exterior regions, as well as the total mass conservation.  In Figure \ref{fig:shrink1}(e), the inner concentration $\psi^{\rm in}$ is approaching the equilibrium value $\psi_{\star}^{\rm in}=0.1$, while outer concentration $\psi^{\rm out}$ roughly stays at the equilibrium value $\psi_{\star}^{\rm out}=0.8$.

\begin{figure}
\includegraphics[width=\textwidth]{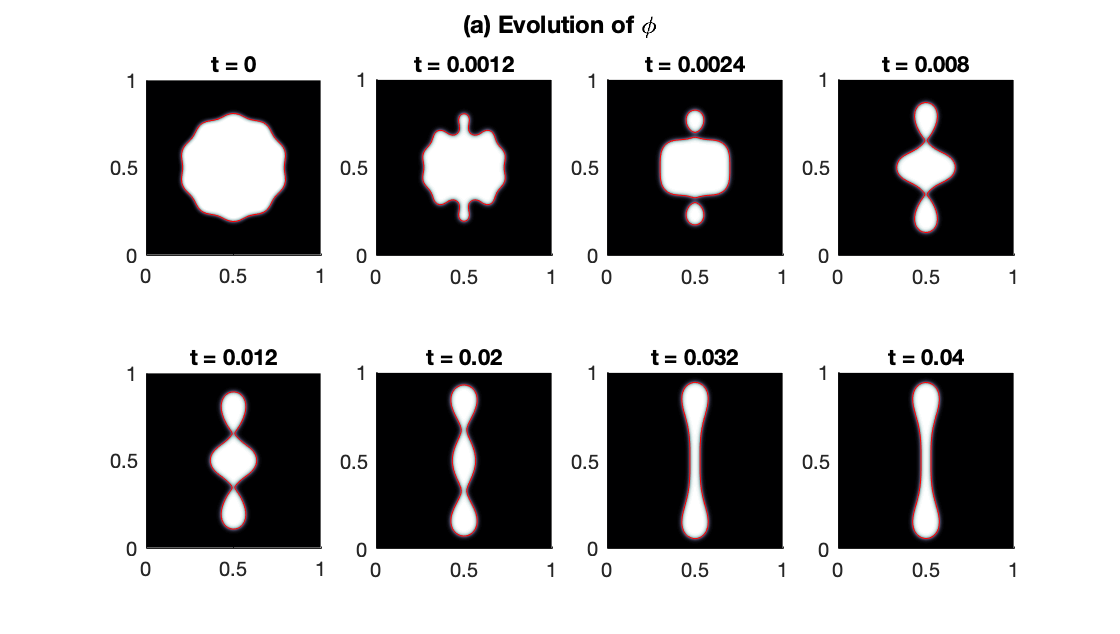}\\
\vspace{-0.2cm} 
\includegraphics[width=\textwidth]{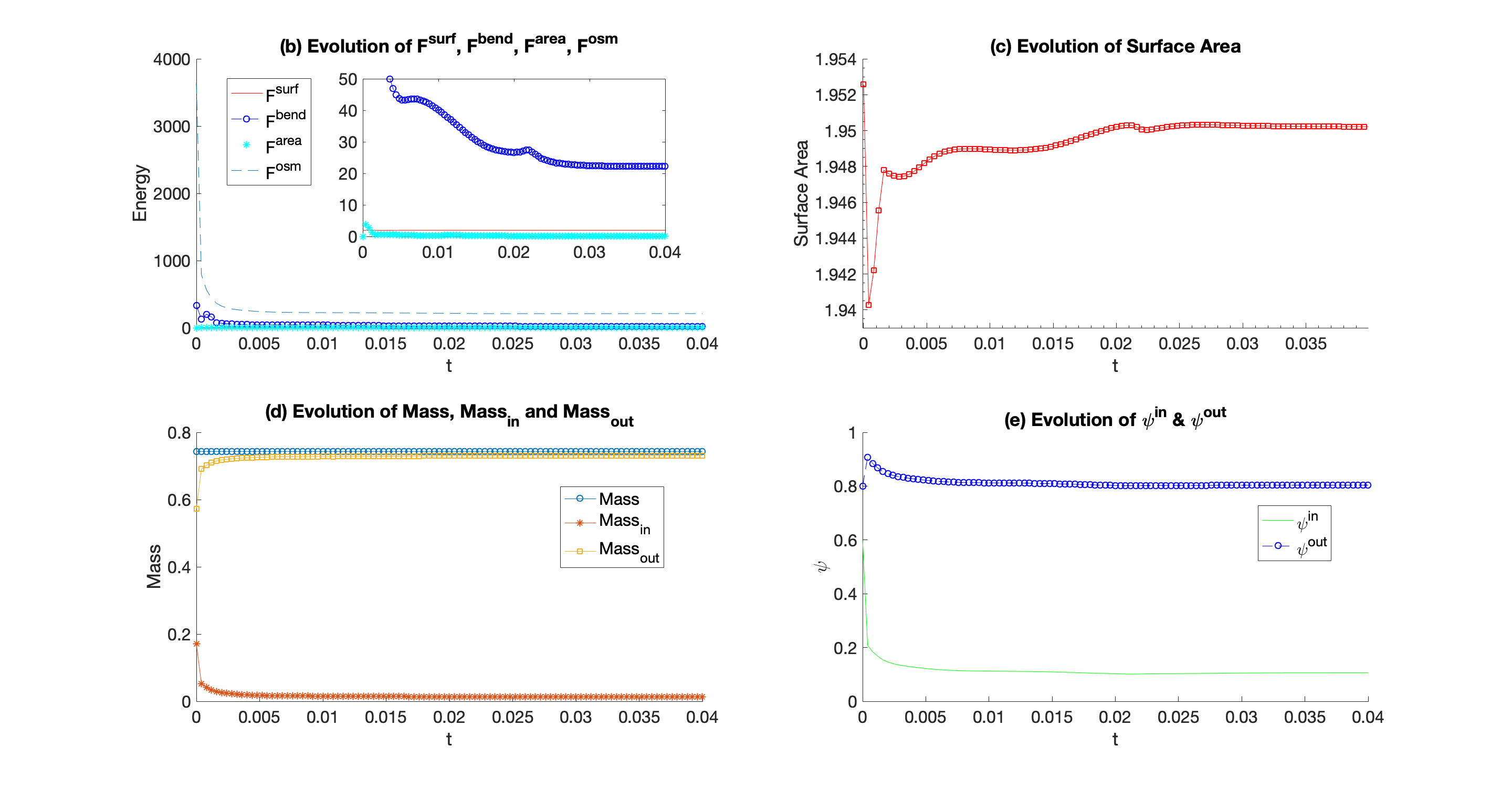}
\caption{The evolution of $\phi$, energy, surface area, mass, and concentration in another shrinkage case. The initial condition and the parameters are the same as those in Figure \ref{fig:shrink1} except $\gamma_{\rm bend}= 1$. (a) shows the similar shrinkage effect as in Figure \ref{fig:shrink1}: the interface become wrinkled with permanent area, but the shape is different due to a much greater $F^{\text{bend}}$. In (b), decline of the osmotic energy causes the inner region shrinking. $F^{\text{bend}}$ decreases, $F^{\text{surf}}$ remains roughly a constant and $F^{\text{area}}$ stays close to 0. (c)-(d) present the interface area constraint and mass conservation. In (e), $\psi^{\rm in}$ approaches the equilibrium value $\psi_{\star}^{\rm in}=0.1$, while $\psi^{\rm out}$ stays at the equilibrium value $\psi_{\star}^{\rm out}=0.8$.}
\label{fig:shrink2}
\end{figure}

Next, we study the shape evolution when $\gamma_{\rm bend}=1$, which is a ten times bending
energy compared with the previous computation. In Figure \ref{fig:shrink2},  we find that the inner phase shrinks with a very different pattern due to a much larger bending energy. 
The energy changes are similar to the ones in Figure \ref{fig:shrink1}(b). In Figure \ref{fig:shrink2}(c)-(d), we present the surface area constraint and mass conservation. In Figure \ref{fig:shrink2}(e), we show that the concentrations inside and outside the interface approach the equilibrium values $\psi_{\star}^{\rm in}$ and $\psi_{\star}^{\rm out}$, respectively. 

\begin{figure}[H]
\includegraphics[width=0.8\textwidth]{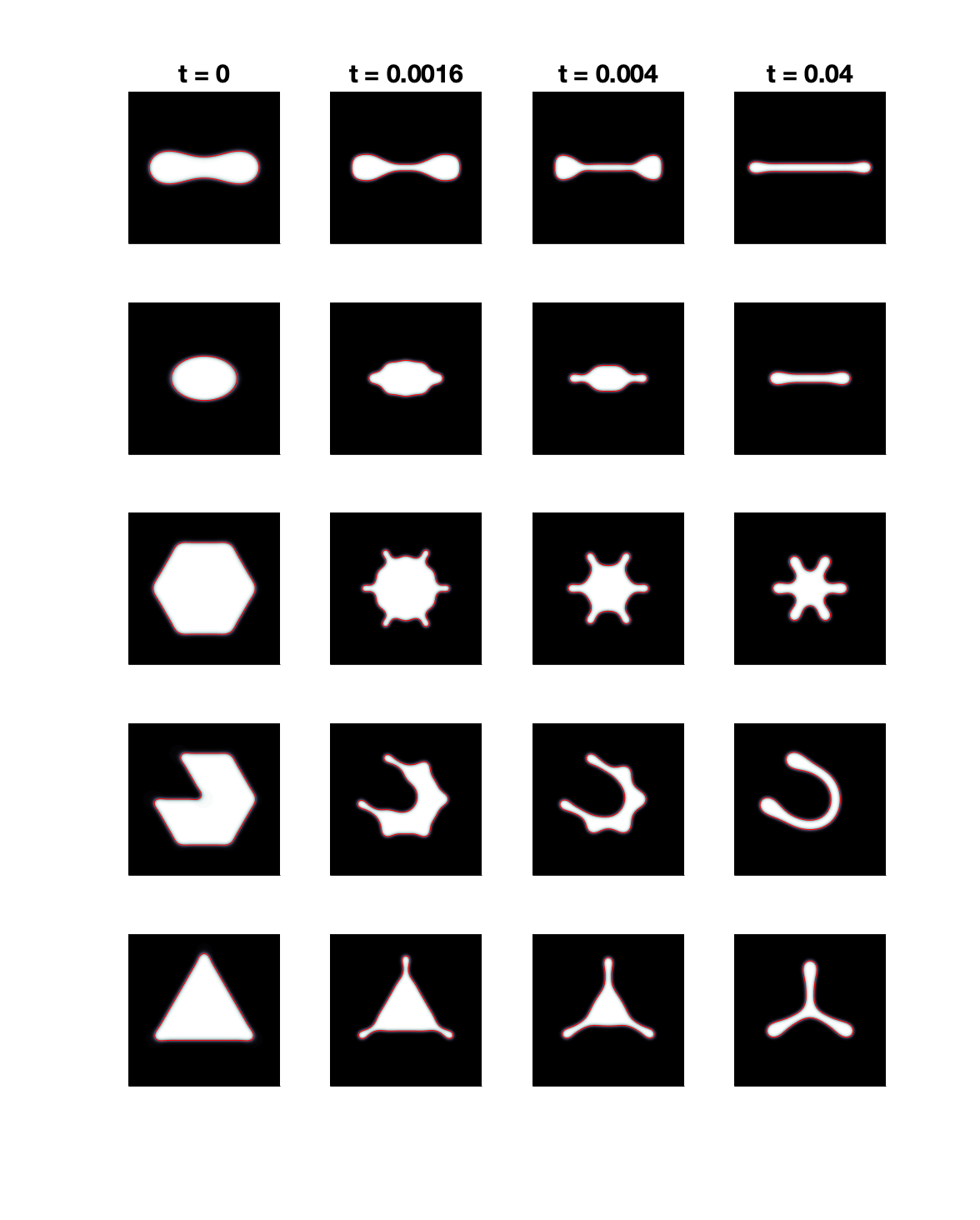}
\centering
\caption{Shape evolution of other shrinkage examples with different initial morphology but same parameters as those used in Figure \ref{fig:shrink1}. 
They have a similar pattern: the inner region shrinks for the decline of osmotic energy and the sharp corners shrink to form finger-like structure; the shrinking process is with preserved surface area and mass; the concentrations of the inner and outer phases evolve to the corresponding equilibrium values.}
\label{fig:shrink_other}
\end{figure}

In Figure \ref{fig:shrink_other}, we present shape evolution of several other shrinkage examples with different initial morphology.  They all follow the similar pattern: the white region will shrink for the decline of osmotic energy; the sharp corners will shrink to form finger-like structures; the shrinking process is always associated  with arclength constraint and total mass conservation. Note that the concentrations of the inner and outer phases evolve to their corresponding equilibrium values.

\section{Conclusion}
In this paper, we have developed a phase field model for vesicle growth or shrinkage based on osmotic pressure that arises due to a chemical potential gradient. The primary contribution is that we defined a neat form of osmotic energy \eqref{eq:fosm} with two simple quadratic functions \eqref{f_inout}. Another novelty is the control conditions derived by the common tangent construction for growth or shrinkage in Section \ref{sec:3}. We directly discretized the model (a coupled  Allen-Cahn and Cahn-Hilliard equations), and implemented a nonlinear FAS multigrid method for computation. We showed detailed convergence tests  and presented numerically the evolution of growth and shrinkage cases. It turns out that the proposed model provide an effective way for describing vesicle growth or shrinkage in different situations by changing the parameters in the multiple kinds of free energies. 

Though we have only done the numerical computing without stability and solvability analysis in this paper, we plan to construct a new scheme based on a convex splitting of the discrete energy, by which both the unconditional unique solvability and energy stability of the numerical scheme are assured \cite{Wise2010, Hu2009}.
In addition, we plan to extend this work to 3D and a more complex evolution by adding a Stokes-like equation, following some ideas used in  \cite{Yang2015,Yang2021}. This fluid-structure interaction type model will enable a full realization of vesicle dynamics, though challenging due to the two extra variables: velocity and pressure.

\section{acknowledgements}
SL acknowledges the support from the National Science Foundation (NSF), Division of Mathematical Sciences grant DMS-1720420. SL was also partially supported by grant ECCS-1307625. SW acknowledges support from NSF grant, DMS 2012634. JL acknowledges partial support from the NSF through grants DMS-1714973, DMS-1719960, and DMS-1763272 and the Simons Foundation (594598QN) for a NSF-Simons Center for Multiscale Cell Fate Research. JL also thanks the National Institutes of Health for partial support through grants 1U54CA217378-01A1 for a National Center in Cancer Systems Biology at UC Irvine and P30CA062203 for the Chao Family Comprehensive Cancer Center at UC Irvine.

\bibliographystyle{ieeetr}
\bibliography{bibfile}

\begin{thebibliography}{10}

\bibitem{Elani2015}
Y.~Elani, R.~V. Law, and O.~Ces, ``Vesicle-based artificial cells as chemical
  microreactors with spatially segregated reaction pathways,'' {\em Nature
  communications}, 2015.

\bibitem{Albert2002}
B.~Alberts, A.~Johnson, J.~Lewis, M.~Raff, K.~Roberts, and P.~Walter, {\em
  Molecular Biology of the Cell. 4th edition.}
\newblock New York: Garland Science, 2002.

\bibitem{Strange2004}
K.~Strange, ``Cellular volume homeostasis,'' {\em Advances in Physiology
  Education}, vol.~28, no.~4, pp.~155--159, 2004.
\newblock PMID: 15545344.

\bibitem{BAUMGARTEN2012}
C.~M. Baumgarten and J.~J. Feher, ``Chapter 16 - osmosis and regulation of cell
  volume,'' in {\em Cell Physiology Source Book (Fourth Edition)}
  (N.~Sperelakis, ed.), pp.~261--301, San Diego: Academic Press, fourth
  edition~ed., 2012.

\bibitem{wiki}
{LadyofHats}, ``Tonicity --- {W}ikipedia{,} the free encyclopedia,'' 17
  February 2007.

\bibitem{VEERAPANENI2009A}
S.~K. Veerapaneni, D.~Gueyffier, D.~Zorin, and G.~Biros, ``A boundary integral
  method for simulating the dynamics of inextensible vesicles suspended in a
  viscous fluid in 2d,'' {\em Journal of Computational Physics}, vol.~228,
  no.~7, pp.~2334--2353, 2009.

\bibitem{VEERAPANENI2009B}
S.~K. Veerapaneni, D.~Gueyffier, G.~Biros, and D.~Zorin, ``A numerical method
  for simulating the dynamics of 3d axisymmetric vesicles suspended in viscous
  flows,'' {\em Journal of Computational Physics}, vol.~228, no.~19,
  pp.~7233--7249, 2009.

\bibitem{SOHN2010119}
J.~S. Sohn, Y.-H. Tseng, S.~Li, A.~Voigt, and J.~S. Lowengrub, ``Dynamics of
  multicomponent vesicles in a viscous fluid,'' {\em Journal of Computational
  Physics}, vol.~229, no.~1, pp.~119--144, 2010.

\bibitem{Salac2011}
D.~Salac and M.~Miksis, ``A level set projection model of lipid vesicles in
  general flows,'' {\em Journal of Computational Physics}, vol.~230, no.~22,
  pp.~8192--8215, 2011.

\bibitem{Sohn2012}
J.~Sohn, S.~Li, X.~Li, and J.~Lowengrub, ``Axisymmetric multicomponent
  vesicles: A comparison of hydrodynamic and geometric models,'' {\em
  International Journal for Numerical Methods in Biomedical Engineering},
  vol.~28, pp.~346 -- 368, 03 2012.

\bibitem{Shuwang2012}
S.~Li, J.~Lowengrub, and A.~Voigt, ``Locomotion, wrinkling, and budding of a
  multicomponent vesicle in viscous fluids,'' {\em Communications in
  Mathematical Sciences}, vol.~10, 06 2012.

\bibitem{FRANK2013}
F.~Hau{\ss}er, W.~Marth, S.~Li, J.~Lowengrub, A.~Rätz, and A.~Voigt,
  ``Thermodynamically consistent models for two-component vesicles,'' {\em
  International Journal of Biomathematics and Biostatistics}, vol.~2,
  pp.~19--48, January-June 2013.

\bibitem{Kai2016}
K.~Liu, C.~Hamilton, J.~Allard, J.~Lowengrub, and S.~Li, ``Wrinkling dynamics
  of fluctuating vesicles in time-dependent viscous flow,'' {\em Soft Matter},
  vol.~12, pp.~5663--5675, 2016.

\bibitem{Salac2018}
P.~Gera and D.~Salac, ``Three-dimensional multicomponent vesicles: dynamics and
  influence of material properties,'' {\em Soft Matter}, vol.~14,
  pp.~7690--7705, 2018.

\bibitem{Salac2022}
P.~Gera, D.~Salac, and S.~E. Spagnolie, ``Swinging and tumbling of
  multicomponent vesicles in flow,'' {\em Journal of Fluid Mechanics},
  vol.~935, p.~A39, 2022.

\bibitem{Kai2014}
K.~Liu and S.~Li, ``Nonlinear simulations of vesicle wrinkling,'' {\em
  Mathematical Methods in The Applied Sciences}, vol.~37, pp.~1093--1112, 2014.

\bibitem{Kai2017}
K.~Liu, G.~R. Marple, J.~Allard, S.~Li, S.~Veerapaneni, and J.~Lowengrub,
  ``Dynamics of a multicomponent vesicle in shear flow,'' {\em Soft Matter},
  vol.~13, pp.~3521--3531, 2017.

\bibitem{Layton2006}
A.~Layton, ``Modeling water transport across elastic boundaries using an
  explicit jump method,'' {\em SIAM Journal on Scientific Computing}, vol.~28,
  p.~2189, 2006.

\bibitem{Vogl2014}
C.~Vogl, M.~Miksis, S.~Davis, and D.~Salac, ``The effect of glass-forming
  sugars on vesicle morphology and water distribution during drying,'' {\em
  Journal of The Royal Society Interface}, vol.~11, no.~99, 2014.

\bibitem{Jayathilake20101}
P.~Jayathilake, Z.~Tan, B.~Khoo, and N.~Wijeysundera, ``Deformation and osmotic
  swelling of an elastic membrane capsule in stokes flows by the immersed
  interface method,'' {\em Chemical Engineering Science}, vol.~65, no.~3,
  p.~1237–1252, 2010.

\bibitem{Jayathilake20102}
P.~Jayathilake, B.~Khoo, and Z.~Tan, ``Effect of membrane permeability on
  capsule substrate adhesion: Computation using immersed interface method,''
  {\em Chemical Engineering Science}, vol.~65, no.~11, p.~3567–3578, 2010.

\bibitem{Mori2011}
Y.~Mori, C.~Liu, and R.~Eisenberg, ``A model of electrodiffusion and osmotic
  water flow and its energetic structure,'' {\em Physica D: Nonlinear
  Phenomena}, vol.~240, p.~1835–1852, 2011.

\bibitem{Yao2017}
L.~Yao and Y.~Mori, ``A numerical method for osmotic water flow and solute
  diffusion with deformable membrane boundaries in two spatial dimension,''
  {\em J. Comput. Phys.}, pp.~728--746, 2017.

\bibitem{WANG2020}
X.~Wang, X.~Gong, K.~Sugiyama, S.~Takagi, and H.~Huang, ``An immersed boundary
  method for mass transfer through porous biomembranes under large
  deformations,'' {\em Journal of Computational Physics}, vol.~413, p.~109444,
  2020.

\bibitem{Quaife2021}
B.~Quaife, A.~Gannon, and Y.-N. Young, ``Hydrodynamics of a semipermeable
  inextensible membrane under flow and confinement,'' {\em Phys. Rev. Fluids},
  vol.~6, p.~073601, Jul 2021.

\bibitem{PESKIN1977220}
C.~S. Peskin, ``Numerical analysis of blood flow in the heart,'' {\em Journal
  of Computational Physics}, vol.~25, no.~3, pp.~220--252, 1977.

\bibitem{DU2004450}
Q.~Du, C.~Liu, and X.~Wang, ``A phase field approach in the numerical study of
  the elastic bending energy for vesicle membranes,'' {\em Journal of
  Computational Physics}, vol.~198, no.~2, pp.~450--468, 2004.

\bibitem{Xiaoqiang2008}
X.~Wang and Q.~Du, ``Modelling and simulations of multi-component lipid
  membranes and open membranes via diffuse interface approaches,'' {\em Journal
  of mathematical biology}, vol.~56, pp.~347--71, 04 2008.

\bibitem{Lowengrub2009}
J.~Lowengrub, A.~Rätz, and A.~Voigt, ``Phase-field modeling of the dynamics of
  multicomponent vesicles: Spinodal decomposition, coarsening, budding, and
  fission,'' {\em Physical review. E, Statistical, nonlinear, and soft matter
  physics}, vol.~79, p.~031926, 04 2009.

\bibitem{Gu2016ATP}
R.~Gu, X.~Wang, and M.~D. Gunzburger, ``A two phase field model for tracking
  vesicle–vesicle adhesion,'' {\em Journal of Mathematical Biology}, vol.~73,
  pp.~1293--1319, 2016.

\bibitem{Giga2017}
M.~Giga, A.~Kirshtein, and C.~Liu, {\em Variational Modeling and Complex
  Fluids. In: Giga Y., Novotny A. (eds) Handbook of Mathematical Analysis in
  Mechanics of Viscous Fluids}.
\newblock Springer, 2017.

\bibitem{Kobayashi2010}
R.~Kobayashi, ``A brief introduction to phase field method,'' {\em AIP
  Conference Proceedings}, vol.~1270, pp.~282--291, July 2010.

\bibitem{Shen2012}
J.~Shen, X.~Yang, and q.~Wang, ``Mass and volume conservation in phase field
  models for binary fluids,'' {\em Communications in Computational Physics},
  vol.~13, 01 2012.

\bibitem{Provatas2010}
N.~Provatas and K.~Elder, ``Phase‐field methods in materials science and
  engineering,'' {\em Wiley-VCH}, 2010.

\bibitem{Bartels2015}
S.~Bartels, {\em The Allen–Cahn Equation. In: Numerical Methods for Nonlinear
  Partial Differential Equations. Springer Series in Computational
  Mathematics}.
\newblock Springer, Cham, 2015.

\bibitem{Lee2014}
D.~Lee, J.-Y. Huh, D.~Jeong, J.~Shin, A.~Yun, and J.~Kim, ``Physical,
  mathematical, and numerical derivations of the cahn–hilliard equation,''
  {\em Computational Materials Science}, vol.~81, pp.~216--225, 2014.

\bibitem{Wise2010}
S.~Wise, ``Unconditionally stable finite difference, nonlinear multigrid
  simulation of the cahn-hilliard-hele-shaw system of equations,'' {\em J. Sci.
  Comput.}, vol.~44, pp.~36--68, 2010.

\bibitem{Trottenberg2005}
U.~Trottenberg, C.~Oosterlee, and A.~Schüller, {\em Multigrid}.
\newblock Academic Press, 2005.

\bibitem{Henson2003}
V.~Henson, ``Multigrid methods for nonlinear problems: An overview,''
  vol.~5016, 12 2002.

\bibitem{Kay2006}
D.~Kay and R.~Welford, ``A multigrid finite element solver for the
  cahn-hilliard equation,'' {\em J. Comput. Phys.}, vol.~212, p.~288–304,
  2006.

\bibitem{Cahn1958}
J.~Cahn and J.~Hilliard, ``Free energy of a nonuniform system. i. interfacial
  free energy,'' {\em J. Chem. Phys.}, vol.~28, no.~2, p.~258–267, 1958.

\bibitem{Du2005}
Q.~Du, C.~Liu, R.~Ryham, and X.~Wang, ``A phase field formulation of the
  willmore problem,'' {\em Nonlinearity}, vol.~18, no.~3, pp.~1249--1267, 2005.

\bibitem{Pelton2019}
A.~D. Pelton, {\em Phase Diagrams and Thermodynamic Modeling of Solutions}.
\newblock Elsevier, 2019.

\bibitem{Hu2009}
Z.~Hu, S.~Wise, C.~Wang, and J.~Lowengrub, ``Stable and efficient
  finite-difference nonlinear-multigrid schemes for the phase-field crystal
  equation,'' {\em J. Comput. Phys.}, vol.~228, p.~5323–5339, 2009.

\bibitem{Yang2015}
R.~Chen, G.~Ji, X.~Yang, and H.~Zhang, ``Decoupled energy stable schemes for
  phase-field vesicle membrane model,'' {\em Journal of Computational Physics},
  vol.~302, pp.~509--523, 2015.

\bibitem{Yang2021}
X.~Yang, ``Numerical approximations of the navier–stokes equation coupled
  with volume-conserved multi-phase-field vesicles system: Fully-decoupled,
  linear, unconditionally energy stable and second-order time-accurate
  numerical scheme,'' {\em Computer Methods in Applied Mechanics and
  Engineering}, vol.~375, p.~113600, 2021.

\end{thebibliography}

	\end{document}